\documentclass[aop,MSNbibl,dvips]{arximspdf}
\usepackage{graphics}
\doi{10.1214/10-AOP532}
\volume{38}
\issue{5}
\pubyear{2010}
\firstpage{1901}
\lastpage{1923}

\makeatletter
\newtheorem{theorem}{Theorem}
\newtheorem{lem}{Lemma}
\newtheorem{prop}{Proposition}
\newtheorem{cor}{Corollary}
\newproclaim{defn}{Definition}
\newproclaim{exmp}{Example}
\newproclaim{rem}{Remark}

\makeatother

\begin{document}
\begin{frontmatter}

\title{Heat kernel estimates for the fractional Laplacian~with
Dirichlet conditions}
\runtitle{Heat kernel estimates}

\begin{aug}
\author[a]{\fnms{Krzysztof} \snm{Bogdan}\corref{}\thanksref{t1}\ead[label=em1]{krzysztof.bogdan@pwr.wroc.pl}},
\author[a]{\fnms{Tomasz} \snm{Grzywny}\thanksref{t2}\ead[label=em2]{tomasz.grzywny@pwr.wroc.pl}}
\and
\author[a]{\fnms{Micha\l{}} \snm{Ryznar}\thanksref{t2}\ead[label=em3]{michal.ryznar@pwr.wroc.pl}}
\runauthor{K. Bogdan, T. Grzywny and M. Ryznar}
\thankstext{t1}{Supported in part by Grant MNiSW N N201 397137.}
\thankstext{t2}{Supported in part by Grant MNiSW N N201 373136.}
\affiliation{Wroc\l{}aw University of Technology}
\address[a]{K. Bogdan\\ T. Grzywny\\ M. Ryznar\\ Institute of Mathematics
and\\\quad Computer Sciences\\
Wroc\l{}aw University of Technology\\ ul. Wybrze\.{z}e Wyspia\'
{n}skiego 27\\ 50-370 Wroc\l{}aw\\Poland\\
\printead{em1}\\\phantom{E-mail: }\printead*{em2}\\ \phantom
{E-mail: }\printead*{em3}} %adresu isvedimo komanda gale!
\end{aug}

% HISTORY:
\received{\smonth{7} \syear{2009}}
\revised{\smonth{1} \syear{2010}}

% ABSTRACT
%
\begin{abstract}
We give sharp estimates for the heat kernel of the fractional
Laplacian with Dirichlet
condition for a general class of
domains including Lipschitz domains.
\end{abstract}

% KEYWORDS
%


\end{frontmatter}

%s1 ###
\section{Introduction}\label{sec1}

Explicit {sharp} estimates for the Green function of
the Laplacian in $C^{1,1}$ domains were completed in 1986 by Zhao \cite
{MR842803}.
% (see also \cite{MR657523, MR902422}).
Sharp estimates of the Green function of Lipschitz domains
were given in 2000 by Bogdan \cite{MR1741527}.
Explicit {qualitatively} sharp estimates for the classical heat
kernel in $C^{1,1}$ domains were established in 2002 by Zhang
\cite{MR1900329}.
% (see also \cite{MR1974093, MR879702},
%and \cite{MR2255354, MR1968386} for further extensions).
Qualitatively sharp heat kernel estimates in Lipschitz domains
were given in 2003 by Varopulous \cite{MR1969798}.
The development of the boundary potential theory of the fractional
Laplacian follows a parallel path. Green function estimates
were obtained in 1997 and 1998 by Kulczycki
\cite{MR1490808} and Chen and Song \cite{MR1654824} for $C^{1,1}$ domains,
% (see \cite[Corollary 1.8]{CKS2008} for the case of dimension one, see
%also \cite{MR1825645}),
and in 2002 by Jakubowski for Lipschitz domains \cite{MR1991120}.
%(see also \cite{MR2213639, MR2182071}).
In 2008 Chen, Kim and Song \cite{CKS2008} gave
sharp explicit estimates for the heat kernel $p_D(t,x,y)$ of the
fractional Laplacian
on $C^{1,1}$ domains $D$.
The main contribution of the present paper is the following
%approximate
result.
%factorization.
% of $p^D(t,x,y)$.
% of in geometrically regular $D$.
%, a factorization which verifies for the cones and $C^{1,1}^$ domains.
%, suggested by the case of the
%right circular cones resolved in.
%form of the estimate for a more general class of domains:
%, see (\ref{eq:CKS}) below.
%the right circular cones.

\begin{theorem}\label{theorem:oppz}
If $D$ is $\kappa$-fat, then there is $C=C(\alpha,D)$ such that
% constants $, c_2$ exist such that for all $x, y\in\Rd$,
%
%e1 ###
\begin{equation}
\label{eq:get}
C^{-1} P^{x}(\tau_{D}>t)P^{y}(\tau_{D}>t)\le\frac{p_{D}(t, x,
y)}{p(t,x,y)}
\le C P^{x}(\tau_{D}>t)
P^{y}(\tau_{D}>t)
\end{equation}
for $0<t\le1$ and $x,y\in D$.
\end{theorem}

Here $p(t,x,y)$ is the heat kernel of the fractional Laplacian on
%the whole space
${\mathbb{R}^{d}}$, and
\begin{eqnarray*}
\label{eq:sp}
P^x(\tau_D>t)=\int_{\mathbb{R}^{d}}p_D(t,x,y)\,{d}y
\end{eqnarray*}
defines the {\it survival
probability} of the corresponding isotropic $\alpha$-stable L\'evy
process in $D$.
%We note that t
The result applies also to unbounded domains, in particular, to domains
above the graph of a Lipschitz function, where we can take {\it arbitrary}
$t>0$. In fact, (\ref{eq:get}) holds with $C=C(\alpha,d,\kappa)$ under
the mere condition that $D$ is $(\kappa,t^{1/\alpha})$-fat at $x$ and
at $y$; see
% the definitions and the proof of Theorem~\ref{theorem:oppz}
Sections~\ref{sec:es} and \ref{sec:sld} for definitions and results.
For {\it exterior} domains we have a
result free from local geometric assumptions:

\begin{cor}\label{theorem:ext}
If $\operatorname{diam}( D^c)<\infty$, then {\rm(\ref{eq:get})} holds with
$C=C(\alpha,d)$
for all $t> \operatorname{diam}(D^c)^\alpha$ and $x,y \in D$.
\end{cor}

For exterior domains of class $C^{1,1}$ a more explicit estimate
is given in Theorem~\ref{cor:c2} below.
We also like to note that a useful variant of Theorem~\ref{theorem:oppz}
is given in
Theorem~\ref{theorem:bLd}.

Expression (\ref{eq:get}) is motivated by these applications of
the semigroup property of~$p_D$:
\[
p_D(2t,x,y)=\int_{\mathbb{R}^{d}}p_D(t,x,z)p_D(t,z,y)\,dz\le
P^x(\tau_D>t)c(t) ,
\]
where $c(t)=\sup_{z,y\in{\mathbb{R}^{d}}}p(t,z,y)\ge \sup
_{z,y\in{\mathbb{R}^{d}}}p_D(t,z,y)$
[see (\ref{eq:gg})], and
% Analogously,
%
\begin{eqnarray*}p_D(3t,x,y)&=&\int\int
p_D(t,x,z)p_D(t,z,w)p_D(t,w,y)\,dw\,dz
\\
& \le&
P^x(\tau_D>t)c(t)P^y(\tau_D>t) .
\end{eqnarray*}
The latter inequality is quite satisfactory for $x=y$, because $c(t)=p(t,x,x)$.
{Off-diagonal} $(x,y)$ in (\ref{eq:get})
require, however, a deeper analysis. Our proof of (\ref{eq:get}) is
based on the
% following properties of the isotropic stable L\'evy processs:
boundary Harnack principle (BHP) \cite{MR2365478} (see also earlier
\cite{MR1719233}), a version of the Ikeda--Watanabe \cite{MR0142153} formula
(\ref{eq:IWf}),
 scaling (\ref{eq:scD}) and comparability of $p$
with its L\'evy measure (\ref{eq:lm}); see (\ref{eq:cpn}). Counterparts
of these are important in view of possible generalizations.

In what follows (\ref{eq:get}) and analogous {\it sharp estimates}
will be written as
\[
p_{D}(t, x, y)\stackrel{C}{\approx}
P^{x}(\tau_{D}>t) p(t,x,y) P^{y}(\tau_{D}>t) ,
\]
meaning that either ratio of the sides is bounded by a number $C\in
(0,\infty)$, and~$C$ does not depend on the variables shown (here: $t$, $x$, $y$). We
will skip $C$ from notation if unimportant for our goals.
% We
%generally fix the capitalized constants $C, C_1, C_2,\ldots$
%throughout the text, but the lower case constants $c, c_1, c_2,\ldots$
%may change value from place to place. Unless stated otherwise,
%constants depend {\it only} on $d$, $\alpha$ and $\kappa$. This will
%sometimes be emphasized by writing, e.g., $C=C(d,\alpha,\kappa)$.
%Notation $C(a,b,\ldots,c)$ means that $C$ may so be chosen to depend
%only on $a,b,\ldots,c$.

Let $\delta_D(x)=\operatorname{dist}(x,D^c)$. As mentioned above, domains $D$
of class $C^{1,1}$ enjoy
the following sharp and explicit estimate of
Chen, Kim and Song \cite{CKS2008}:
%
%e2 ###
\begin{equation}
\label{eq:CKS}
\, \frac{p_D(t,x,y)}{p(t,x,y)}\approx
\biggl(1\wedge\frac{\delta_D^{\alpha/2}(x)}{t^{1/2}}\biggr)
\biggl(1\wedge\frac{\delta_D^{\alpha/2}(y)}{t^{1/2}}\biggr)
 ,\qquad 0<t\le
1 , x,y\in{\mathbb{R}^{d}} .
\end{equation}
We note that (\ref{eq:CKS}) agrees with (\ref{eq:get}) because by
\cite{2009KBTGcm}, Corollary 1,
\begin{eqnarray*}
\label{eq:1}
P^x(\tau_D>t)\approx
1\wedge\frac{\delta_D^{\alpha/2}(x)}{t^{1/2}}
\qquad\mbox{for } 0<t\le1 ,  x,y\in{\mathbb{R}^{d}} .
\end{eqnarray*}
In fact, starting with (\ref{eq:get}), we are able to recover and
strengthen (\ref{eq:CKS}), with a simpler proof; see Example~\ref{ex:c11}
and Proposition~\ref{cor:espc11} below. We note that
(\ref{eq:get}) was conjectured in \cite{2009KBTGcm} based on the cases
of $C^{1,1}$ domains \cite{CKS2008} and circular cones
\cite{2009KBTGcm}. We should also mention that the Gaussian estimates
of Varopoulous \cite{MR1969798} have a shape similar to (\ref{eq:get}),
in particular, they involve the survival probability.
Thus, the present paper builds on the evidence accumulated in
\cite{CKS2008,MR1969798} and \cite{2009KBTGcm}.
We also note that the {upper} bound in (\ref{eq:CKS}) was
proved in 2006 by Siudeja for semibounded convex domains (\cite{MR2255353}, Theorem~1.6), and stated for general convex domains
in \cite{MR2255353}, Remark~1.7. Some of our present techniques were inspired
%It appears that an impulse for our
%by the proof of (\ref{eq:CKS})
%was given
by \cite{MR2231884}, Theorem~4.2, of Kulczycki and Siudeja,
\cite{MR2438694}, Proposition 2.9, of Ba\~nuelos and
Kulczycki, and \cite{MR2075671}, Section 4, of Bogdan and Ba\~nuelos.
%A similar but weaker upper bound was earlier given in
%, see also \cite{MR2280260, 2008-MK, MR0208681}.
%Generally, the subject is far from exhausted--and it seems manageable
%with the existing techniques.

%We complement (\ref{eq:get}) with estimates of the survival
%probability for $\kappa$-fat domains $D$.
It is a consequence of Lemma~\ref{lem:etd1r} below
%that $x\mapsto P^x(\tau_D>1)$ is comparable with harmonic functions at
%the boundary of $D$, and
that we can apply BHP
\cite{MR2365478,MR1719233} to conveniently estimate $P^x(\tau_D>1)$
by some kernel
functions of $D$, namely, by the Martin
kernel with the pole at infinity {\it or} the expected survival time
[we use {\it scaling} to estimate $P^x(\tau_D>t)$ for general $t>0$].
The estimate and the resulting bounds for the heat kernel are collected
in Theorem~\ref{theorem:bLd}, followed by a number of applications.
%For special Lipschitz domains sharp estimates of the survival
%probability can be expressed in terms
%of the Martin kernel, see Theorem~\ref{theorem:oppsld} below.
%We give applications
%to
%the ball, the complement of the ball etc.
In particular, we
%As a consequence
give a simple proof of the main result of \cite{2009KBTGcm}
for the circular cones~$V$:
%
%e3 ###
\begin{equation}
\label{eq:etdG2}
\frac{p_{V}(t,x,y)}{p(t,x,y)}
\approx
\frac{(
1\wedge{\delta_V(x)/t^{1/\alpha}}
)^{\alpha/2}}{
(
1\wedge{|x|/t^{1/\alpha}}
)^{\alpha/2-\beta}}\frac{
(
1\wedge{\delta_V(y)/t^{1/\alpha}}
)^{\alpha/2}}{
(
1\wedge{|y|/t^{1/\alpha}}
)^{\alpha/2-\beta}} .
\end{equation}
Here $\beta\in[0,\alpha)$ is a characteristic of the cone, and
{\it all}\/ $t>0$ and $x,y\in{\mathbb{R}^{d}}$ are allowed.
%Bounded $\kappa$-fat domains are resolved in Theorem~\ref{theorem:bLd}
%below.
We should add to (\ref{eq:get}), (\ref{eq:CKS}) and (\ref{eq:etdG2})
that \cite{MR0119247,MR2013738}
%
%e4 ###
\begin{equation}
\label{eq:setd}
p_t(x) \stackrel{c}{\approx}
\frac{t}{|x|^{d+\alpha}} \land t^{-d/\alpha} ,\qquad t>0 , x\in
{\mathbb{R}^{d}} .
\end{equation}
Here $c=c(\alpha,d)$, meaning that $c\in(0,\infty)$ may be so chosen
to depend only on $d$ and $\alpha$.
%Noteworthy, all the above-mentioned estimates for bounded $C^{1,1}$
%domains have the same form as for the ball
% (in this connection compare \cite[Corollary 1.2]{CKS2008} with
%We also like to note that the right circular cones are merely special
%Lipschitz domains,
%but a number of techniques and explicit formulas make them
%an interesting and important test case (see \cite{MR1465162, MR863716,
%MR1062058, MR0474525, MR799436}).
%We hope to encourage a further study of Lipschitz and more general
%domains for
%stable and other jump-type processes \cite{MR1942325, MR2238879,
% MR2085428, MR2006232}.
We like to note that the estimates for general $\kappa$-fat
%(in particular: for Lipschitz)
domains cannot be as explicit as those for $C^{1,1}$
domains. In particular, the decay rate $\beta$
at the vertex of a cone delicately depends on the
aperture of the cone; see
\cite{MR2075671,2009KBTGcm,MR2213639} (see also \cite{MR1741527}).
Nevertheless, Lipschitz domains offer a natural
setting for studying the boundary behavior of the Green function and
the heat kernel for both the Brownian motion and the isotropic $\alpha
$-stable L\'evy
processes.
This is due to the {scaling}, the rich range of
asymptotic behaviors depending on the local geometry of the domain's
boundary, connections to
the boundary Harnack principle and approximate factorization of the Green
function, and applications in the perturbation
theory of generators, in particular, via the 3G Theorem
\cite{MR2075671,MR1741527,MR1671973,MR2207878,MR842803} and 3P
Theorem \cite{MR2283957}.
The $\kappa$-fat sets are a convenient generalization of Lipschitz
domains, with similar features.
It is noteworthy that (\ref{eq:get}) is an approximate factorization
of the
heat kernel (see also \cite{MR1741527,MR2365478} in this connection).

We should add that the $C^{1,1}$ condition specifies the geometry of a
domain only in bounded scales (see Definition~\ref{defn:3}).
This renders the range of time in (\ref{eq:CKS}) restricted to
$0<t\le1$.
%Cones are also examples of unbounded domains, which are only partially
%resolved by
%the results of \cite{CKS2008, CKS2008censored} (note that (
%valid only for bounded times).
In what follows we will also study the probability of survival for
large times (and unbounded domains). This is straightforward for
special Lipschitz domains (thus for circular cones),
but less so for general $\kappa$-fat or $C^{1,1}$ domains.
As an interesting case study we consider
domains with bounded
complement (i.e., exterior domains) of class $C^{1,1}$.
These have distinctive geometries at infinity and at the
boundary, resulting in nontrivial completion of (\ref{eq:CKS}).
%especially for $\alpha\ge 1=d$ (the recurrent case).
%Sharp explicit estimates of the survival probability for exterior
%domains are
%given in Section~\ref{sec:ex}. The estimates complement
%Theorem~\ref{theorem:ext}. We note that general $C^{1,1}$ domains are also
%studied in
%Section~\ref{sec:ex}.
We remark that exterior $C^{1,1}$
domains in dimension $d>\alpha$ have been recently studied in
\cite{2009ZCJT}, too.
We also remark that \cite{MR2386098}, Theorem~4.4 bounds
the survival probability of the relativistic process in a half-line,
and \cite{2009KKMS} gives an explicit formula
for the transition density of the killed Cauchy process ($\alpha=1$)
on the half-line.
Regarding other recent estimates
\cite{MR2492992,MR2320691,MR2357678,MR2430977,MR2238934}
for the transition density and potential kernel of jump-type processes,
we need to point out that generally these only concern processes
without killing.
Killing corresponds to the Dirichlet ``boundary'' condition (analogous
to the negative Schr\"odinger perturbation
\cite{MR1825645,2008KBTJ}) and it severely influences the
asymptotics of the transition density and Green function.
Needless to say, the asymptotics are crucial for solving the
Dirichlet problem \cite{MR2417435,MR2386098}.
%In fact, the Dirichlet problem for more general nonlocal operators,
%here
%exemplified by the fractional Laplacian, can be conveniently studied
%within a similar probabilistic setting.
%As we shall see, the heat kernel of the fractional Laplacian in the
%right
%circular cones has a
%power-type asymptotics at infinity, and it decays like the distance to
%the
%boundary to the power $\alpha/2$ except at the vertex, where it decays
%with the rate of $\beta\in(0,\alpha)$.

We like to mention possible applications and further directions of
research. The estimate (\ref{eq:get}) fits well into the technique of
Schr\"odinger perturbations of \cite{2008KBTJ}, which should
produce straightforward consequences. Also, the distribution of $\tau
_D$, given by (\ref{eq:IWf}) below, can be estimated by using (\ref{eq:get}).
Further, we conjecture that for certain domains $D$,
$\lim p_D(t,x,y)/P^x(\tau_D>t)$ exists as $x$ approaches a boundary
point of $D$.
This may lead to representation theorems for nonnegative parabolic
functions of the fractional Laplacian (compare \cite{MR2365478}, Theorems~2 and~3) and construction of excursion laws. We need to remark
here that our estimates are inconclusive about the (irregular \cite
{MR2365478}) boundary points of $D$, but we conjecture that (\ref
{eq:get}) indeed extends to $\partial D$.
Finally, it seems important to understand the
behavior of $p_D(t,x,y)$ for domains which are rather small at a
boundary point or at infinity.
%For instance, if $D$ is
%bounded then $p_D(t,x,y)$ decays exponentially when $t\to\infty$, due
%to intrinsic ultracontractivity
In this connection we refer the interested reader to the recent study
of {\it intrinsic
ultracontractivity} by Kwa\'snicki \cite{2009-MK-pa}; see also
\cite{2009KBTGcm,CKS2008,MR1643611} and the notion of
{\it inaccessibility} in \cite{MR2365478}.

Our general references to the {boundary potential theory} of the fractional
Laplacian are \cite{MR1671973} and \cite{MR2365478}. We also refer
the reader to \cite{KB-TB-MR-TK-RS-ZV} for a broad non-technical
overview of the methods and goals of the theory.

The paper is composed as follows. In Section~\ref{sec:prel} we recall
basic facts about the killed isotropic $\alpha$-stable
L\'evy processes.
In Section~\ref{sec:es} we prove Theorem~\ref{theorem:oppz} and
Corollary~\ref{theorem:ext}.
In Section~\ref{sec:sld} we state and prove Theorem~\ref{theorem:bLd}
and give applications to specific domains.
In particular, we strengthen (\ref{eq:CKS}) and part of the results of
\cite{CKS2008} (see Proposition~\ref{cor:espc11}, Theorem~\ref{cor:c2}
and Corollary~\ref{cor:c4}), and we discuss exterior $C^{1,1}$ domains
in dimension $d=1<\alpha$.

%s2 ###
\section{Preliminaries}\label{sec:prel}

In what follows, ${\mathbb{R}^{d}}$ denotes the Euclidean space of
dimension $d\ge
1$, $dy$ is the Lebesgue measure on ${\mathbb{R}^{d}}$, and $0<\alpha<2$.
Our primary analytic data are as follows: a nonempty open set $D\subset
{\mathbb{R}^{d}}$ and the L\'evy measure given by density
function
%
%e5 ###
\begin{equation}
\label{eq:lm}
\nu(y)=\frac{2^{\alpha}\Gamma((d+\alpha)/2)}{\pi^{d/2}|\Gamma
(-\alpha/2)|}
|y|^{-d-\alpha} .
\end{equation}
The coefficient in (\ref{eq:lm}) is such that
%
%e6 ###
\begin{equation}
\label{eq:trf}
\int_{\mathbb{R}^{d}}[1-\cos(\xi\cdot y)]\nu
(y)\,dy=|\xi|^\alpha
,\qquad
\xi\in{\mathbb{R}^{d}} .
\end{equation}
For (smooth compactly supported) $\phi\in C^\infty_c({\mathbb{R}^{d}})$,
the fractional Laplacian is
\begin{eqnarray*}
\Delta^{\alpha/2}\phi(x) &=&
\lim_{\varepsilon\downarrow0}\int_{|y|>\varepsilon}
[\phi(x+y)-\phi(x)]\nu(y)\,dy ,
\qquad
x\in{\mathbb{R}^{d}}
\end{eqnarray*}
(see \cite{MR1671973,KB-TB-MR-TK-RS-ZV} for a broader setup).
%, i.e. $\phi:\Rd\to\R$ is smooth and compactly supported on $\Rd$.
If $r>0$ and $\phi_r(x)=\phi(r x)$, then
%
%e7 ###
\begin{equation}
\label{eq:scfr}
\Delta^{\alpha/2}\phi_r(x)=r^\alpha\Delta^{\alpha/2}\phi(rx)
,\qquad
x\in{\mathbb{R}^{d}} .
\end{equation}
We let $p_t$ be the smooth real-valued function on ${\mathbb{R}^{d}}$
with Fourier transform,
%
%e8 ###
\begin{equation}
\label{eq:dpt}
\int_ {\mathbb{R}^{d}}p_t(x)e^{ix\cdot\xi} \,dx=e^{-t|\xi|^\alpha
} ,\qquad t>0
, \xi\in
{\mathbb{R}^{d}} .
\end{equation}
In particular, the maximum of $p_t$ is
$p_t(0)=2^{1-\alpha}\pi^{-d/2}\alpha^{-1}\Gamma(d/\alpha)/\Gamma
(d/2)t^{-d/\alpha}$.
According to (\ref{eq:trf}) and the L\'evy--Khinchine formula,
$\{p_t\}$ is a probabilistic convolution semigroup with L\'evy measure
$\nu(y)\,dy$; see \cite{MR2013738,MR1739520} or
\cite{KB-TB-MR-TK-RS-ZV}. We have the following scaling property,
%
%e9 ###
\begin{equation}
\label{eq:sca}
p_t(x)=t^{-d/\alpha}p_1(t^{-1/ \alpha}x) ,\qquad t>0 , x\in
{\mathbb{R}^{d}} ,
\end{equation}
which may be considered a consequence of (\ref{eq:dpt}).
It is noteworthy that by (\ref{eq:setd}) we have
%
%e10 ###
\begin{equation}\label{eq:ppp}
p_t(x)\approx p_{2t}(x) ,\qquad t>0 ,  x\in{\mathbb{R}^{d}} .
\end{equation}
%
% The semigroup $P_t f(x)=\int_ \Rd f(y) p_t(y-x) dy$ has
% $\Delta^{\alpha/2}$ as infinitesimal generator (\cite{MR0481057},
% \cite{MR1336382}, \cite{MR1671973}, \cite{MR1873235}). Put
% differently,
We denote
\[
p(t,x,y)=p_t(y-x) ,
\]
and we have
%
%e11 ###
\begin{equation}
\label{eq:fsol}
\int_{s}^\infty\int_{ {\mathbb{R}^{d}}}
p(u-s,x,z)[
\partial_u\phi(u,z)+\Delta^{\alpha/2}_z \phi(u,z)] \,dz\,du
= -\phi(s,x) ,
\end{equation}
where $s\in\mathbb{R}$, $x\in{\mathbb{R}^{d}}$, and $\phi\in
C^\infty_c(\mathbb{R}\times
{\mathbb{R}^{d}})$;
see, for example, \cite{2008KBTJ}, (36).

We define the isotropic $\alpha$-stable L\'evy process $(X_t,P^x)$ by
stipulating transition probability
\[
P_t(x,A)=\int_A p(t, x,y)\,dy ,\qquad t>0 , x\in{\mathbb{R}^{d}}
,  A\subset
{\mathbb{R}^{d}} ,
\]
initial distribution $P^x(X(0)=x)=1$, and c\'adl\'ag paths. Thus,
$P^x$, $E^x$ denote the distribution and expectation for the process
starting at $x$.
% It is well-known that $(X_{t},P^{x})$ is strong Markov with respect
% to the so-called standard filtration \cite{MR1406564, MR0264757}.
We define the {\it time of the first exit} from $D$, or {\it survival time},
\[
\tau_D=\inf\{t>0\dvtx   X_t\notin D\} ,
\]
and the {\it time of first hitting} $D$,
\[
T_D=\inf\{t>0\dvtx X_t\in D\} .
\]
We define, as usual,
\[
p_D(t,x,y)=p(t,x,y)-E^x[\tau_D<t;  p(t-\tau_D, X_{\tau_D},y)],\qquad t>0
, x,y\in{\mathbb{R}^{d}} .
\]
%
% see, e.g., \cite{CKS2008, MR1671973}.
We have that
%
%e12 ###
\begin{equation}\label{eq:gg}
0\le p_D(t,x,y)=p_D(t,y,x)\le p(t,x,y) ,
\end{equation}
hence,
%
%e13 ###
\begin{equation}
\label{eq:9.5}
\int p_D(t,x,y)\,dy=\int p_D(t,x,y)\,dx\le1 .
\end{equation}
If $x\in D^c$ is {\it regular}\/ for the Dirichlet problem on $D$ \cite
{MR2365478}, that is, $P^x(\tau_D=0)=1$, then $p_D(t,x,y)=0$ and (\ref
{eq:get}) is trivially satisfied. By this remark, if all the points of
$\partial D$ are regular for $D$, then we can write $x,y\in{\mathbb
{R}^{d}}$ in
Theorem~\ref{theorem:oppz}, instead of~\mbox{$x,y\in D$}.
The remark also applies to Examples \ref{ex:1}--\ref{ex:8} in
Section~\ref{sec:sld}.
By the strong Markov property, $p_D$ is the transition density of the
isotropic stable process {\it killed} when leaving $D$, meaning that we have
the following Chapman--Kolmogorov equation,
\[
\int_{\mathbb{R}^{d}}p_D(s,x,z)p_D(t,z,y)\,dz=p_D(s+t,x,y) ,\qquad
s,t>0 ,
x,y\in{\mathbb{R}^{d}} ,
\]
and for nonnegative or bounded (Borel) functions $f \dvtx {\mathbb
{R}^{d}}\to\mathbb{R}$,
\[
\int_{\mathbb{R}^{d}}f(y)p_D(t,x,y)\,dy=E^x[\tau_D<t;  f(X_t)] ,
\qquad t>0 ,  x\in{\mathbb{R}^{d}} .
\]
For $s\in\mathbb{R}$, $x\in{\mathbb{R}^{d}}$, and $\phi\in
C^\infty_c(\mathbb{R}\times D)$, we
have
\begin{eqnarray*}
\label{eq:fsolD}
\int_{s}^\infty\int_{D}
p_D(u-s,x,z)[
\partial_u\phi(u,z)+\Delta^{\alpha/2}_z \phi(u,z)]\, dz\,du
= -\phi(s,x) ,
\end{eqnarray*}
which extends (\ref{eq:fsol}) and justifies calling $p_D$ the heat
kernel of the (Dirichlet) fractional Laplacian {\it on} $D$. It is
well known that $p_D$ is jointly continuous and positive for
$(t,x,y)\in(0,\infty)\times D\times D$. We have a scaling property,
$p_{rD}(r^\alpha t, rx,ry)=r^d p_D(t,x,y)$, $r>0$, or
%
%e14 ###
\begin{equation}
\label{eq:scD}
p_D(t,x,y)=t^{-d/\alpha}
p_{t^{-1/\alpha}D}(1,t^{-1/\alpha}x,t^{-1/\alpha}y) ,\qquad t>0 ,
x,y\in{\mathbb{R}^{d}} ,
\end{equation}
in agreement with (\ref{eq:sca}) and (\ref{eq:scfr}). Thus,
$P^{rx}(\tau_{rD}>r^\alpha t)=P^x(\tau_D>t)$, or
%
%e15 ###
\begin{equation}
\label{eq:scDs}
P^x(\tau_D>t)=\int_{{\mathbb{R}^{d}}}
p_D(t,x,y)\,dy=P^{t^{-1/\alpha}x}(\tau_{t^{-1/\alpha}D}>1) .
\end{equation}

\begin{rem}\label{rem:cnu}
For $c>0$ consider $\tilde{\nu}=c\nu$, the corresponding heat
kernels $\tilde{p}$, $\tilde{p}_D$, probability and expectation
$\tilde{P}^x$, $\tilde{E}^x$. Clearly,
$\tilde{p}_D(t,x,y)={p}_D(ct,x,y)$.
\end{rem}

The Green function of $D$ is defined as
%
%e16 ###
\begin{equation}
\label{eq:12.5}
G_D(x,y)=\int_0^{\infty} p_D(t,x,y)\,dt ,
\end{equation}
and scaling of $p_D$ yields the following scaling of $G_D$,
%
%e17 ###
\begin{equation}
\label{eq:3}
G_{rD}(rx,ry)=r^{\alpha-d}G_D(x,y) .
\end{equation}

A result of Ikeda and Watanabe \cite{MR0142153} asserts that for $x\in
D$ the $P^x$-distribution of $(\tau_D, X_{{\tau_D}-},X_{\tau_D})$
restricted to $X_{{\tau_D}-}\neq X_{\tau_D}$ is given by the density
function
%
%e18 ###
\begin{equation}
\label{eq:IWf}
(s,u,z)\mapsto p_D(s,x,u)\nu(z-u) .
\end{equation}
For geometrically nice domains, for example, for the ball,
$P^x(X_{{\tau_D}-}\neq X_{\tau_D})=1$ for $x\in D$ \cite{MR2365478}, and
then
% The harmonic measure $\omega^x_D(\cdot)$ at $x$ of $D$ is the
% $P^x$-distribution of $X_{\tau_D}$. If $P^x(X_{{\tau_D}-}\neq
% X_{\tau_D})=0$,
by (\ref{eq:12.5}) and (\ref{eq:IWf}) the $P^x$-distribution of
$X_{\tau_D}$ has the density function given by the Poisson kernel,
%
%e19 ###
\begin{equation}
\label{eq:Pk}
P_D(x,z)=\int_D G_D(x,u)\nu(z-u)\,du .
\end{equation}
%
% $P_D(x,z)$ is called the Poisson kernel of $D$.
For $x_0\in{\mathbb{R}^{d}}$ and $r>0$ we consider the ball
$B(x_0,r)=\{x\in{\mathbb{R}^{d}}
\dvtx
|x-x_0|<r\}$ and $B^c(x_0,r)=\{x\in{\mathbb{R}^{d}}\dvtx |x-x_0|>r\}$
(open complement
of a ball).

There is a constant $C$ depending only on $d$,
$\alpha$ and $p$, such that
%
%e20 ###
\begin{equation}\label{eq:bhp}
P_U(x_1, y_1)P_U(x_2, y_2) \stackrel{C}{\approx}
P_U(x_1, y_2)P_U(x_2, y_1)  ,
\end{equation}
whenever $U\subset B(x_0,r)\subset{\mathbb{R}^{d}}$ is open, $0<
p<1$, $r>0$,
$x_0\in{\mathbb{R}^{d}}$,
$x_1, x_2 \in U \cap B(x_0,rp)$, and $y_1, y_2\in B(x_0,r)^c$.
This boundary Harnack principle (BHP) follows from
\cite{MR2365478}, Lemma~7 and the proof of Theorem~1, and it is
essentially an approximate
factorization of $P_U$. We encourage the interested reader to
directly verify the estimate in the special case of
(\ref{eq:poisson:ball}) below.
% We claim that (\ref{eq:bhp}) remains true for all $x_1, x_2 \in
%B(x_0,r/2)$, and
% $y_1, y_2\in B(x_0,r)^c$. Indeed, if $y_j\in D$ for $j=1$ or $2$,
% and $x_i$ is irregular for (the Dirichlet problem on) $D$
% (\cite{MR2365478}) for $i=1$ and $2$, then $P_D(x_1,y_j)=\infty$
%because
% $G_D(x_i,u)>0$ for $u\in D$, and the singularity of $\nu$ is
% non-integrable. This makes (\ref{eq:bhp}) trivial in this case.
% If $i=1$ or $2$ and $x_i$ is regular for $D$, then
% $G_D(x_i,u)=0$ in (\ref{eq:bhp}), and
% $P_D(x_1,y_j)=0$ for $j=1$ and $2$, thus proving (\ref{eq:bhp}).
% In particular, if $d=1<\alpha$ then each point of $D^c$ is regular for
% $D$ and we are done. For $\alpha<d$ we consider
% $D_\varepsilon=D\cup B(x_1,\varepsilon)\cup B(x_2,\varepsilon)$, and
% let $\varepsilon\to0^+$. We have that $G_{D_\varepsilon}\downarrow
% G_D$, and by monotone convergence $P_{D_\varepsilon}\downarrow P_D$.
% The constant $C$ in (\ref{eq:bhp}) is independent of the domain, so
% the comparability is preserved in the limit.

The Green function and Poisson kernel of $B(x_0,r)$
are known explicitly:
%
%e22 ###
%e21 ###
\begin{eqnarray}\label{wfg}
G_{B(x_0,r)}(x,v)&=&
\mathcal{B}_{d,\alpha}
|x-v|^{\alpha-d}\int_{0}^{w}\frac{s^{\alpha/2-1}}
{(s+1)^{d/2}} \,ds ,
\\ \label{eq:poisson:ball}
P_{B(x_0,r)}(x, y) &=&
\mathcal{C}_{d, \alpha}
\biggl[\frac{r^2 - |x-x_0|^2}{|y-x_0|^2 -
r^2}\biggr]^{\alpha/ 2} \frac{1}{|x - y|^d}  ,
\end{eqnarray}
where $\mathcal{B}_{d,\alpha}=\Gamma(d/2)/(2^{\alpha}\pi^{d/2}[\Gamma
(\alpha/2)]^{2})$,
$\mathcal{C}_{d,\alpha}=\Gamma(d/2)\pi^{-1-d/2}\sin(\pi\alpha/2)$,
\[
w=(r^2-|x-x_0|^{2})(r^2-|v-x_0|^{2})/|x-v|^{2} ,
\]
$|x-x_0|<r$, $|v-x_0|<r$, and $|y-x_0|\ge  r$;
see \cite{MR0126885,bibRm}.
Thus,
%
%e23 ###
\begin{equation}
\label{eq:tpk}
P^x\bigl(\bigl|X_{\tau_{B(0,1)}}\bigr|>R\bigr)=\int_{|y|\ge
R}P_{B(0,1)}(x,y)\,dy\approx
\frac{(1-|x|)^{\alpha/2}}
{R^{\alpha}} ,
\end{equation}
where $x\in B(0,1)$ and $ R\ge 2$. Also, for $|x-x_0|\le r$
we have
\cite{MR1825645}
%
%e24 ###
\begin{equation}
\label{eq:sBr}
E^x \tau_{B(x_0,r)}(x)=\frac{2^{1-\alpha}\Gamma(d/2)}
{\alpha\Gamma((d+\alpha)/2)\Gamma(\alpha
/2)}(r^2-|x-x_0|^2)^{\alpha/2} .
\end{equation}

All the sets and
functions considered below are Borelian.
{\it Positive} means {\it strictly positive}.
{\it Domain} means a nonempty open set (connectedness need not be
assumed in this theory).
% for special Lipschitz domains, and sharp {\it explicit} estimates
% for the survival probability for exterior domains.

%s3 ###
\section{Factorization}\label{sec:es}

We consider nonempty open set $D\subset{\mathbb{R}^{d}}$.
% and $\overline{D}=D\cup\partial D$.
%
\begin{defn}\label{def:UB}
Let $x\in D$, $r>0$ and $0<\kappa\le1$. We say that $D$ is
$(\kappa,r)$-fat at $x$ if there is a ball $B(A,\kappa r)\subset
D\cap B(x,r)$. If this is true for every $x\in{D}$, then we say
that $D$ is $(\kappa,r)$-fat. We say that $D$ is $\kappa$-fat if
there is $R>0$ such that $D$ is $(\kappa,r)$-fat for all $r\in
(0,R]$.
\end{defn}

\begin{rem}\label{rem:ball}
The ball is $1/2$-fat.
\end{rem}

\begin{defn}\label{def:U}
Given $B(A,\kappa)\subset D\cap B(x,1)$, we consider $U= D\cap
B(x,|x-A|+\kappa/3)$, $B_1=B(A,\kappa/3)\subset U$ and
$B_2=B(A',\kappa/6)$ such that $B(A',\kappa/3)\subset
B(A,\kappa)\setminus U$; see the picture:
\end{defn}

\begin{figure}[h]

\includegraphics{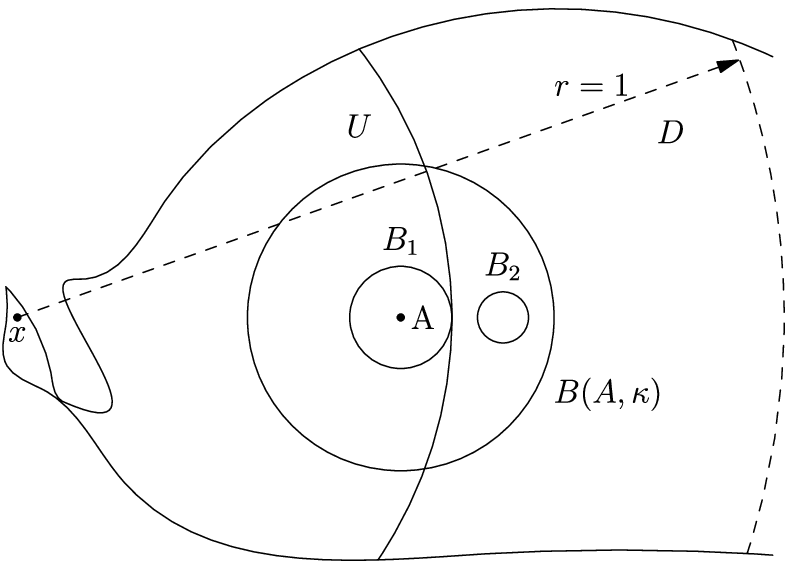}

  %\caption{}
\end{figure}

\begin{lem}\label{lem:etd1r}
There is $C=C(\alpha,d,\kappa)$ such that if $D$ is
$(\kappa,1)$-fat at $x$, then
%
%e25 ###
\begin{equation}\label{eq:etd1r1}
P^x(\tau_{D}>1/3)\le C
P^x(\tau_{D}>3) .
\end{equation}
\end{lem}

\begin{pf}
% In what follows $x\in U_1$.
Consider $x\in D$ and $B(A,\kappa)$ and $U$ as above. For
$|x-A|<\kappa/2$,
\begin{eqnarray*}
1&\ge & P^x(\tau_D>1/3)\ge  P^x(\tau
_D>3)
\\
&\ge &
P^x\bigl(\tau_{B(x,\kappa/2)}>3\bigr)=P^0\bigl(\tau_{B(0,\kappa/2)}>3\bigr)>0 ,
\end{eqnarray*}
and
(\ref{eq:etd1r1}) is proved. We will now assume that $|x-A|\ge
\kappa/2$. We note that
%
%e26 ###
\begin{equation}
\label{eq:rsp}
P^x(\tau_{D}>1/3) \le P^x(\tau_{U}>1/3)+P^x(X_{\tau_{U}}\in
D) .
\end{equation}
We have $P^x(X_{\tau_{U}}\in D)=\int_D P_U(x,y)\,dy$. Indeed,
if $B=B(x,|x-A|+\kappa/3)$ as in the definition of $U$, then
$P^x(X_{\tau_{U}}\in\partial U\cap D)\le
P^x(X_{\tau_{B}}\in\partial B)=0$; see the discussion preceding
(\ref{eq:Pk}) above. Similarly, $P^x(X_{\tau_{U}}\in B_2)$ is an
integral of the Poisson kernel $P_U$. We consider BHP
for
$x_1=x$, $x_2=A$, $p=1-\kappa/3>(1-\kappa)/(1-\kappa+\kappa/3)$.
Integrating (\ref{eq:bhp}) on $D$ and $B_2$, we obtain
\[
\frac{P^x(X_{\tau_{U}}\in D)}
{P^{A}(X_{\tau_{U}}\in D)} \le
c\frac{P^x(X_{\tau_{U}}\in B_2)}{P^{A}(X_{\tau_{U}}\in B_2)} .
\]
We note that (the denominator) $P^{A}(X_{\tau_{U}}\in B_2)\ge
P^{A}(X_{\tau_{B_1}}\in B_2)\ge  c>0$ [see (\ref{eq:poisson:ball})],
therefore, $P^x(X_{\tau_{U}}\in D)\le c P^x(X_{\tau_{U}}\in
B_2)$. We
also observe that $u \mapsto\int_{B_2}\nu(y-u)\,dy$ is bounded away
from zero and infinity on $U$.
% The process does not hit $\partial B(Q,r_3)$ when leaving neither
% $B(Q,r_3)$ nor, consequently, $U_3$.
By (\ref{eq:Pk}),
\begin{eqnarray*}
\label{eq:cesp}
P^x(X_{\tau_{U}}\in B_2)= \int_{U}G_{U}(x,u)\int_{B_2}\nu(y-u)\,dy\,du
\approx\int_{U}G_{U}(x,u)\,du =E^x\tau_{U} .
\end{eqnarray*}
Clearly, $P^x(\tau_{U}>1/3)\le3 E^x\tau_{U}$. By (\ref
{eq:rsp}), $
P^x(\tau_{D}>1/3)\le c E^x\tau_{U} $. By the strong Markov property,
\begin{eqnarray*}
E^x\tau_{U}
&\le& cP^x(X_{\tau_{U}}\in B_2)
\le c E^x\bigl[ X_{\tau_{U}}\in B_2;
P^{X_{\tau_{U}}}\bigl(\tau_{B(X_{\tau_{U}},\kappa/6)}>3\bigr)\bigr]
\\
&\le& c P^x(\tau_{D}>3) .
\end{eqnarray*}
\upqed\end{pf}

\begin{rem}\label{comp}
If $D$ is $(\kappa,1)$-fat at $x$, then by the above proof we have
%
%e27 ###
\begin{equation}
\label{eq:cw}
\hspace*{8pt}P^x(\tau_{D}>1/3) \approx
P^x(\tau_{D}>3) \approx
P^x(\tau_{D}>1) \approx
P^x(X_{\tau_{U}}\in D) \approx
E^x\tau_{U}  .
\end{equation}
In fact, we can replace $3$ by any finite ${\mathcal E}\ge 1$, at
the expense of having the comparability between {\it each} pair of expressions
in (\ref{eq:cw}) holding with a constant
$C=C(\alpha,d,\kappa,{\mathcal E })$.
\end{rem}

% The following technical result will be used BASIC LEMMA
%
\begin{lem}\label{lemppu100}
Consider open $D_1, D_3\subset D$ such that ${\rm
dist}(D_1,D_3)>0$. Let $D_2=D\setminus(D_1\cup D_3)$. If $x\in
D_1$ and $y \in D_3$, then
\begin{eqnarray*}
\label{eq:ub}
p_{D}(1, x, y)
\le P^x(X_{\tau_{D_1}}\in D_2)\sup_{s<1,  z\in D_2} p(s, z, y)
+ E^x \tau_{D_1} \sup_{u\in D_1,  z\in D_3}\nu(z-u)
\end{eqnarray*}
and
\begin{eqnarray*}\label{eq:lb}
p_{D}(1, x, y)\ge  P^x(\tau_{D_1}>1)
 P^y(\tau_{D_3}>1)\inf_{u\in D_1,  z\in D_3}\nu(z-u).
\end{eqnarray*}
\end{lem}

% CHYBA COS TRZEBA ZALOZYC O ZBIORZE $D$: PRZYNAJMNIEJ POWINNO BYC:
%$$p^{D}(t,x,y)=0,\ x\in D^c.$$
%
\begin{pf}
By the strong Markov property,
\[
p_{D}(1, x,
y)=E^x[p_{D}(1-\tau_{D_1}, X_{\tau_{D_1}}, y),\tau_{D_1}<1 ],
\]
which
is
\begin{eqnarray*}
&& E^x[p_{D}(1-\tau_{D_1}, X_{\tau_{D_1}}, y),\tau_{D_1}<1,
X_{\tau_{D_1}}\in D_2 ]
\\
&&\qquad{} +E^x[p_{D}(1-\tau_{D_1}, X_{\tau_{D_1}},
y),\tau_{D_1}<1, X_{\tau_{D_1}}\in D_3 ] = I+\mathit{II} .
\end{eqnarray*}
Clearly,
\[
I\le P^x(X_{\tau_{D_1}}\in D_2)\sup_{s<1,  z\in D_2}
p(s, z, y) .
\]
Consider $D_1$ such that $P^x(X_{\tau_{D_1}}\in\partial D_1\cap
D)=0$, for example, $D_1$ being an intersection of $D$ with a Lipschitz
domain. By (\ref{eq:IWf}), the density function of $(\tau_{D_1},
X_{\tau_{D_1}})$ at $(s, z)$ for $z\in D$ equals
\begin{eqnarray*}
\label{eq:dls}
f^x(s,z)=\int_{{D_1}}p_{{D_1}}(s, x, u)\nu(z-u)\,du .
\end{eqnarray*}
%
% , or L\'evy system,
For $z\in D_3$,
\[
f^x(s,z)=\int_{{D_1}}p_{{D_1}}(s, x, u)\nu(z-u)\,du \le
% \int_{{D_1}}p_{{D_1}}(s, x, u)du \sup_{u\in D_1,  z\in
% D_3}\nu(z-u)=
P^x(\tau_{D_1}>s) \sup_{u\in D_1,  z\in D_3}\nu(z-u)  ,
\]
hence, by (\ref{eq:9.5}),
\begin{eqnarray*}
\mathit{II}&=&\int_0^1\int_{D_3}p_{D}(1-s, z, y)f^x(s,z)\,dz\, ds
\\
&\le& \sup_{u\in D_1,  z\in D_3}\nu(z-u)
\int_0^1\int_{D_3}p_{D}(1-s, z, y)P^x(\tau_{D_1}>s)\,dz\, ds
\\
&\le&
\int_0^1P^x(\tau_{D_1}>s) \,ds \,\sup_{u\in D_1,  z\in D_3}\nu(z-u)
\le
E^x\tau_{D_1}\sup_{u\in D_1,  z\in D_3}\nu(z-u) .
\end{eqnarray*}
The upper bound follows. The case of general $D_1$ follows by
approximating from below, and continuity of $p$ and $\nu$. The lower
bound obtains analogously
\begin{eqnarray*}
\mathit{II}&\ge & \inf_{u\in D_1,  z\in D_3}\nu(z-u)
\int_0^1\int_{D_3}p_{D}(1-s, z, y)P^x(\tau_{D_1}>s)\,dz\, ds
\\
&\ge &
P^x(\tau_{D_1}>1) \inf_{u\in D_1,  z\in D_3}\nu(z-u)
\int_0^1\int_{D_3}p_{D_3}(1-s, z, y)\,dz\, ds .
\end{eqnarray*}
\upqed\end{pf}

\begin{rem}\label{rem:cD1}
Lemma~\ref{lemppu100}
%is fairly general, in particular it
also holds for
$\tilde{\nu}$, $\tilde{p}$, $\tilde{P}^x$ and $\tilde{E}^x$ of
Remark~\ref{rem:cnu}.
\end{rem}

% In applications of Lemma~\ref{lemppu100} we will typically choose
% $D_1$ and $D_3$ such that $\operatorname{dist}(D_3,D_1)>0$ and we will set
% $D_2=D\setminus(D_1\cup D_3)$.

In what follows we will often use the fact that
%
%e28 ###
\begin{equation}
\label{eq:cpn}
1 \wedge\nu(z-u)
\approx
p(1,u,z) .
\end{equation}

\begin{lem}\label{lemppu}
% Assume {\bf A}.
If $D$ is $(\kappa,1)$-fat at $x$ and $y$, then
\begin{eqnarray*}
\label{eq:ppu}
p_{D}(2,x,y)\le C(\alpha, d, \kappa)  P^x(\tau_{D}>2)
P^y(\tau
_{D}>2) p(2,x,y) .
\end{eqnarray*}
\end{lem}

\begin{pf}
% If $\delta_D(x)\ge\kappa/2$ then $P^x(\tau_{D}>2)\ge
% P^0(\tau_{B(0,\kappa/4)}>2)>0$, hence
% \begin{equation}\label{eq:ppun1}
% p_{D}(1,x,y)\le C P^x(\tau_{D}>2) p(1,x,y).
% \end{equation}
If $|x-y|\le8$, then $p(1,x,y)\approx1$, and by the semigroup
property, (\ref{eq:ppp}) and Lemma~\ref{lem:etd1r},
%
%e29 ###
\begin{eqnarray}\label{eq:ppun2}
p_{D}(1,x,y)&=&\int_{\mathbb
{R}^{d}}p_D(1/2,x,z)p_D(1/2,z,y)\,dz\nonumber
\\
&\le&
\sup_{z} p(1/2,z,y) P^x(\tau_{D}>1/2)
\\
&\le& c P^x(\tau_{D}>1)p(1,x,y) .\nonumber
\end{eqnarray}
Here $c=c(\alpha,d,\kappa)$.
% Finally, let $x\in D\cap B(Q, r_1)=U_1$,
If $|x-y|>8$, then we will apply Lemma~\ref{lemppu100} with $D_1=U
=D\cap B(A, |x-A|+\kappa/3)$, as in Definition~\ref{def:U}, and
$D_3= \{z\in D\dvtx |z-x|>|x-y|/2\}$. Since $\sup_{s<1,  z\in D_2}
p(s, z, y)\le c p(1, x, y)$, and $\sup_{u\in D_1,  z\in
D_3}\nu(z-u) \le c p(1, x, y)$ [see (\ref{eq:cpn})], by
Remark~\ref{comp}, we obtain
\begin{eqnarray}\label{eq:ppun3}
 p_{D}(1, x, y)&\le&
cp(1, x, y)[P^x(X_{\tau_{U}}\in D)+ E^x\tau_{U}]\nonumber
\\[-8pt]\\[-8pt]
&\le& cP^x(\tau_{D}>1) p(1,x,y) ,\nonumber
\end{eqnarray}
hence, by (\ref{eq:ppun2}), (\ref{eq:ppun3}), symmetry, the semigroup
property and Lemma~\ref{lem:etd1r},
\begin{eqnarray*}
p_D(2,x,y)&=&\int p_D(1,x,z)p_D(1,z,y)\,dz
\\
&\le&
cP^x(\tau_{D}>1)P^y(\tau_{D}>1)\int p(1,x,z)p(1,z,y)\,dz
\\
&\le&cP^x(\tau_{D}>2)P^y(\tau_{D}>2)p(2,x,y) .
\end{eqnarray*}
%
% \begin{equation}\label{eq:ppun4}
% p_{D}(1,x,y)\le C(\alpha,d,\kappa)
% P^x(\tau_{D}>1/2) p(1,x,y), \quad x,y\in\R.
% \end{equation}
\upqed\end{pf}

Under the assumptions of Lemma~\ref{lemppu},
$\tilde{C}=\tilde{C}(\alpha, d, \kappa)$ exists such that
%
%e30 ###
\begin{equation}
\label{eq:ppu1}
p_{D}(1,x,y)\le\tilde{C}  P^x(\tau_{D}>1) P^x(\tau_{D}>1)
p(1,x,y) .
\end{equation}
Indeed, according to Remark~\ref{rem:cnu}, we consider
$\tilde{\nu}=\frac{1}{2}\nu$ and the corresponding $\tilde{p}$,
$\tilde{p}_D$, $\tilde{P}^x$, obtaining
\begin{eqnarray*}
p_{D}(1,x,y)&=&\tilde{p}_{D}(2,x,y)\le
\tilde{C}\tilde{P}^x(\tau_{D}>2) \tilde{P}^x(\tau_{D}>2)\tilde
{p}(2,x,y)
\\
&=&\tilde{C} {P}^x(\tau_{D}>1) {P}^x(\tau_{D}>1) {p}(1,x,y) .
\end{eqnarray*}
%
% \begin{rem}

% By Remark 2 we have

% \begin{equation}\label{eq:ppu1}
%p_{D}(1,x,y)\le C  P^x(\tau_{D}>2) p(1,x,y),
% \end{equation}

% which shows by applying semigroup property

% $$ p_{D}(2,x,y)\le C  P^x(\tau_{D}>2) P^x(\tau_{D}>2) p(2,x,y)$$

% \end{rem}
%
\begin{lem}\label{lower bound12}
If $r>0$, then there is a constant $C=C(\alpha,d,r)$ such that
\begin{eqnarray*} p_{B(u,r)\cup B(v,r)}(1, u, v)&\ge & C
p(1,u,v)
,\qquad u,v\in{\mathbb{R}^{d}} .
\end{eqnarray*}
\end{lem}

\begin{pf}
For $|u-v|\ge  r/2$ we use (\ref{eq:cpn}) and Lemma \ref{lemppu100}
with $D= B(u,r)\cup B(v,r)$, $D_1= B(u,r/8)$ and $D_3=B(v,r/8)$:
\begin{eqnarray*}p_{B(u,r)\cup B(v,r)}(1, u, v)&\ge &
P^u(\tau_{D_1}>1)P^v(\tau_{D_3}>1)
\inf_{u\in D_1,  z\in D_3}\nu(z-u)
\\
&\ge & c \bigl[P^0\bigl(\tau_{B(0,r/8)}>1\bigr)\bigr]^2p(1,u, v) .
\end{eqnarray*}
For $|u-v|\le r/2$, by (\ref{eq:setd}), we simply have
% we only need to use joint continuity of $p_{B(0,r)}(u, z,v)$
%
\[
p_{B(u,r)\cup B(v,r)}(1, u, v)\ge \inf_{|z|<r/2} p_{B(0,r)}(1,
0, z)\ge  c \ge  c p(1,u,v) .
\]
\upqed\end{pf}

\begin{lem}\label{lemppl}
% Assume {\bf A}.
If $D$ is
$(\kappa,1)$-fat at $x$ and $y$, then
\begin{eqnarray*}
\label{eq:ppl}
p_{D}(3,x,y)\ge  C(\alpha, d, \kappa)  P^x(\tau_{D}>3)
P^y(\tau
_{D}>3) p(3,x,y)  .
\end{eqnarray*}
\end{lem}

\begin{pf}
Consider $U^x$, $B_2^x$, and $U^y$, $B_2^y$, selected
according to Definition~\ref{def:U} for~$x$ and $y$, correspondingly. By
the semigroup property, Lemma~\ref{lower bound12} with $r=\kappa/6$,
and~(\ref{eq:setd}),
\begin{eqnarray*}
p_{D}(3,x,y)&\ge & \int_{B^y_2}\int_{B_2^x}
p_{D}(1,x,u)p_{D}(1,u,v)p_D(1,v,y)\,du\,dv
\\
&\ge &c p(1,x,y)\int_{B_2^x} p_{D}(1,x,u)\,du\,\int
_{B^y_2}p_D(1,v,y)\,dv .
\end{eqnarray*}
For $u\in B_2^x=B(A',\kappa/6)$, by Lemma~\ref{lemppu100} with
$D_1=U^x=U$ and $D_3=B(A',\kappa/4)$, and by Remark~\ref{comp}, we
obtain
\begin{eqnarray*}
p_{D}(1,x,u)&\ge & P^x(\tau_U>1)P^0\bigl(\tau_{B(0,\kappa
/12)}>1\bigr)\inf
_{w\in
U, z\in D_3}\nu(z-w)
\\
&\ge & c P^x(\tau_U>1)\ge  c
P^x(\tau_D>1) .
\end{eqnarray*}
Similarly, $p_D(1,v,y)\ge  c P^y(\tau_D>1)$, hence, by
Lemma~\ref{lem:etd1r}, we have
\begin{eqnarray*}
p_{D}(3,x,y)&\ge & c P^y(\tau_D>1)p(1,x,y)P^x(\tau_D>1)
\\
&\ge & c
P^y(\tau_D>3)p(3,x,y)P^x(\tau_D>3) .
\end{eqnarray*}
\upqed\end{pf}

Under the assumptions of Lemma~\ref{lemppl} we also have that
%
%e31 ###
\begin{equation}\label{eq:ppl1}
p_{D}(1,x,y)\ge \tilde{C}(\alpha, d, \kappa)  P^x(\tau_{D}>1)
P^y(\tau_{D}>1) p(1,x,y)  .
\end{equation}
This is proved analogously to (\ref{eq:ppu1}).

% \begin{rem}
% By Remark 2 we have
% $$ p_{D}(2,x,y)\ge C  P^x(\tau_{D}>2) P^x(\tau_{D}>2) p(2,x,y)$$
% \end{rem}

\begin{pf*}{Proof of Theorem~\ref{theorem:oppz}}
Assume that $R\ge 1$ and $D$ is $(\kappa,r)$-fat for
$0<r\le R$.
If $t^{1/\alpha}\in(0,R]$, then $t^{-1/\alpha}D$ is
$(\kappa,1)$-fat. The estimate (\ref{eq:get}) follows from
(\ref{eq:ppu1}), (\ref{eq:ppl1}) and scaling; see (\ref{eq:scD})
and (\ref{eq:scDs}). In fact, we have $C=C(\alpha,d,\kappa)$ in
(\ref{eq:get}). If $R<1$, then we argue as in the case of
(\ref{eq:ppu1}) $C=C(\alpha,d,\kappa,R)$ or, alternatively, we use
Remark~\ref{rem:inc} below.
\end{pf*}

% In the same way we obtain the following estimate for large $t$.
% \begin{cor}
% If $D$ is $(\kappa,r)$-fat for every $r>R>0$ then (\ref{eq:get})
% holds for all $t>R^\alpha$.
% \end{cor}
%
\begin{pf*}{Proof of Corollary~\ref{theorem:ext}}
Note that $D$ is $(1/4,r)$-fat for $r\ge 2  \operatorname{diam}(D^c)$,
and so we obtain (\ref{eq:get}) for $t\ge 2^{\alpha}  {\rm
diam}(D^c)$ with the same constant $C$. If we consider
$\tilde{\nu}=2^{-\alpha}\nu$ and argue like in the case of
(\ref{eq:ppu1}), then we obtain the wider range of $t$, as in the
statement of Corollary~\ref{theorem:ext}.
\end{pf*}

\begin{rem}
Since the $\kappa$-fatness condition is more restrictive when
$\kappa$ is bigger, the above constants $C=C(\alpha, d, \kappa)$ may be
chosen decreasing with respect to $\kappa$. Also, if $D$ has a {\it
tangent} inner ball of radius $1$ at every boundary point, then the
constants in Lemmas~\ref{lemppu} and \ref{lemppl} depend only
on $\alpha$ and $d$.
\end{rem}

\begin{rem}\label{rem:inc}
If $D$ is $(\kappa,r)$-fat at $x$ and $1\le K<\infty$, then $D$ is
$(\kappa/K,rK)$-fat at $x$. This observation together with scaling
allows to easily {\it increase time}, compare (\ref{eq:ppu1}) or
(\ref{eq:ppl1}), at the expense of enlarging the constants of
comparability. The argument, however, does not allow to decrease
time. Remark~\ref{rem:cnu} is more flexible in this respect.
\end{rem}

% \begin{proof}
% Use scaling property and Lemma \ref{lemppl}.
% \end{proof}

%s4 ###
\section{Applications}\label{sec:sld}

%We consider (arbitrary dimension) $d\ge 1$ and a $(\kappa,r)$ set
%$D\subset\Rd$.

We let $s_D(x)=E_x\tau_D=\int G_D(x,y)\,dv$ if this expectation is
finite for $x\in D$, otherwise we let
$s_D(x)=M_{D}(x)$, the Martin kernel with the pole at infinity for
$D$,
\[
M_D(x)=\lim_{D\ni y, |y|\to\infty}\frac{G_D(x,y)}{G_D(x_0,y)} .
\]
We should note that this (alternative) definition of $s_D$ is natural
in view
of \cite{MR2365478}, Theorem~2. The choice of $x_0\in D$ is merely
a normalization, $M_{D}(x_0)=1$, and will not be reflected in the
notation.
By the scaling of the Green function (\ref{eq:3}), we obtain
%
%e32 ###
\begin{equation}
\label{eq:5}
\frac{s_{rD}(rx)}{s_{rD}(ry)}= \frac{s_{D}(x)}{s_{D}(y)} ,\qquad
x,y\in D, r>0 .
\end{equation}
We denote by $A_{r}(x)$ or $A_r(x,\kappa,D)$ every point $A$ such that
$B(A,\kappa
r)\subset D\cap B(x,r)$,
%If $D$ is $(\kappa,r)$-fat at $x$ then w
as in Definition~\ref{def:UB}.
It is noteworthy that $A_{r}(x)$ approximately dominates $x$ in terms
of the
distance to $\partial D$:
%
%e33 ###
\begin{equation}
\label{eq:6}
\delta_D(A_{r}(x))\approx r\vee\delta_D(x) .
\end{equation}
If $D$ is $(\kappa,1)$-fat at $x$, then $rD$ is $(\kappa,r)$-fat at
$rx$, and (every) $rA_1(x,\kappa,D)$ may serve as $A_r(rx,\kappa,rD)$.
%Such points are often used in the Carleson estimate to bound values
%of nonnegative harmonic functions, say in $U$ of Definition~

\begin{theorem}\label{theorem:bLd}
If $D$ is
($\kappa, t^{1/\alpha}$)-fat at $x$ and $y$, then
%
%e34 ###
\begin{equation}
\label{eq:spkf}
P^x(\tau_D>t)\stackrel{C}{\approx}\frac{s_D(x)}{s_D(A_{t^{1/\alpha
}}(x))} ,
\end{equation}
where $C=C(d, \alpha, \kappa)$ and, furthermore,
%
%e35 ###
\begin{equation}
\label{eq:ppkf}
p_{D}(t, x, y)
\stackrel{C}{\approx}
\frac{s_D(x)}{s_D(A_{t^{1/\alpha}}(x))}
p(t,x,y)
\frac{s_D(y)}{s_D(A_{t^{1/\alpha}}(y))} .
\end{equation}
\end{theorem}

\begin{pf} To verify (\ref{eq:spkf}), we first let $t=1$ and assume
that $D$ is
$(\kappa,1)$-fat at $x$. Let $A=A_1(x)$. If $E^x \tau_D<\infty$,
then we consider
the set $U\subset D$ of Definition~\ref{def:U}, and we obtain
\[
E^x\tau_D=E^x{\tau_U}+E^x s_D(X_{\tau_U}) .
\]
By Remark~\ref{comp},
$E^x{\tau_U}\approx P^x (\tau_D>1)$.
Since $E^A \tau_U\approx1$, we trivially
have
\[
\frac{E^x{\tau_U}}{E^A{\tau_U}}
\approx P^x (\tau_D>1).
\]
Similarly, $P^A(X_{\tau_U}\in D)\approx1$. By BHP and
Remark~\ref{comp}, we obtain
%
%e36 ###
\begin{equation}
\label{eq:2}
\frac{E^x s_D(X_{\tau_U})}{E^A s_D(X_{\tau_U})}
\approx
\frac{P^x(X_{\tau_U}\in D)}{P^A(X_{\tau_U}\in D)}
\approx
P^x(X_{\tau_U}\in D)\approx
P^x(\tau_D>1) .
\end{equation}
%
%This ends the proof of the approximation in (\ref{eq:spkf}) for
%$t=1$.
This yields (\ref{eq:spkf}) in the considered case.
%for $t=1$.
If $E^x \tau_D=\infty$, then $s_D$ is harmonic and we have
$s_D(x)=E^xs_D(X_{\tau_U})$
(see \cite{MR2365478}, Theorem~2 and (77)) and we proceed directly via
(\ref{eq:2}).
The case of general $t$ in (\ref{eq:spkf}) is obtained by the scaling
of (\ref{eq:5}) and~(\ref{eq:scDs}).
Finally,
(\ref{eq:ppkf})
follows from (\ref{eq:spkf}) and Theorem~\ref{theorem:oppz}.
The resulting comparability constants
depend only on
$\alpha$, $d$ and $\kappa$.
\end{pf}
% We should note that $s_D(A_{1}(x))$ is essentially
% constant locally at the boundary of $D$,
% and $1\approx P^x(X_{\tau_U}\in D)\approx s(x)/s(A_{1}(x))$ for
%$x=A_{1}$.
% Thus $s(x)/s(A_{1}(x))$ can be compared to $P^x(X_{\tau_U}\in
% D)$ of Remark~\ref{comp} by using BHP,
% and $$P^x(\tau_D>1)\appc{C} P^x(X_{\tau_U}\in D)\approx
%s(x)/s(A_{1}(x)),\ x\in D.$$
% The constant $C=C(\alpha, d, \kappa)$.

% By the scaling property
% \begin{equation}\label{eq:spc111}
% P^x(\tau_{D}>t)=
% P^{t^{-1/\alpha}x}(\tau_{t^{-1/\alpha}D}>1)=
% \frac{s_{t^{-1/\alpha}D}(t^{-1/\alpha}x)}{s_{t^{-1/\alpha}D}(A^t_{
%  ,\quad x\in D .
% \end{equation}
% where the point $A^t_{\kappa,1}(y)$ is a point $A_{1}(y)$ but for the
%set $t^{-1/\alpha}D$.
% Using the scaling property for $s_{\rho D}(y)= \rho^\alpha s_{ D}(y/

% \begin{equation}\label{eq:spc112}
% P^x(\tau_{D}>t)=
% P^{t^{-1/\alpha}x}(\tau_{t^{-1/\alpha}D}>1)=
% \frac{s_{D}(x)}{s_{D}(t^{1/\alpha}A^t_{\kappa,1}(t^{-1/\alpha}x))}
%  ,\quad x\in D .
% \end{equation}

% Note that $t^{1/\alpha}A^t_{\kappa,1}(t^{-1/\alpha}x)= A_{t^{1/

% \begin{equation}\label{eq:spc113}
% P^x(\tau_{D}>t)\approx
% \frac{s_{D}(x)}{s_{D}(A_{t^{1/\alpha}}(x))}
%  ,\quad x\in D .
% \end{equation}

% Thus, (\ref{eq:spkf}) follows fom Lemma~\ref{lem:etd1r} (Remark~
% and BHP. (\ref{eq:ppkf}) follows from scaling, (\ref{eq:get}) and
% (\ref{eq:spkf}).
% (\ref{eq:ppkflt}) follows from intrinsic ultracontractivity

\begin{rem}\label{rem:ra}
Assume that $D$ is $\kappa$-fat, so that there is $R>0$ such that $D$
is $(\kappa,r)$-fat for every $r\le R$. Then (\ref{eq:spkf}) and
(\ref{eq:ppkf}) hold with $C=C(d,\alpha,\kappa)$ for all $x,y\in D$
and $t\le R^\alpha$.
\end{rem}

Below we give a number of applications.

\begin{exmp}\label{ex:1}
We let $R>0$ and $D=B(0,R)\subset{\mathbb{R}^{d}}$.
By (\ref{eq:sBr}), the expected\vspace*{-1pt} survival time is $s_D(x)\stackrel
{C}{\approx}
\delta_D^{\alpha/2}(x)R^{\alpha/2}$, where $C=C(d,\alpha)$.
%By Remark~\ref{rem:ra} with $\kappa=1/2$ we obtain
By (\ref{eq:6}),
$s_D(A_{t^{1/\alpha}}(x))\stackrel{C}{\approx}
(
t^{1/\alpha}
\vee
\delta_D(x)
) ^{\alpha/2}R^{\alpha/2}$,
therefore, for all $t\le R^{\alpha}$ and $x,y\in{\mathbb{R}^{d}}$,
%
%e37 ###
\begin{equation}\label{eq:ex1_1}
P^x(\tau_{D}>t)\stackrel{C}{\approx}
\frac{\delta_D^{\alpha/2}(x)}{(t^{1/\alpha}\vee\delta_D(x))
^{\alpha/2}}=
\biggl(1 \wedge\frac{\delta_D(x)}{t^{1/\alpha}}\biggr)^{\alpha/2}
\end{equation}
and
%
%e38 ###
\begin{equation}
\label{eq:ex1_2}
p_{D}(t, x, y)
\stackrel{C}{\approx}
\biggl(1 \wedge\frac{\delta_D^{\alpha/2}(x)}{ t^{1/2}}\biggr)
p(t,x,y)
\biggl(1 \wedge\frac{\delta_D^{\alpha/2}(y)}{t^{1/2}}\biggr) .
\end{equation}
\end{exmp}

To be explicit, $\delta_{B(0,R)}(x)=(R-|x|)\vee0$, and $\delta
_{B(0,R)^c}(x)=(|x|-R)\vee0$, and (\ref{eq:ex1_1}), (\ref{eq:ex1_2})
on $D^c$ follow because all $x\in D^c$ are regular for $D$.

\begin{exmp}\label{tail_point}
Let $D\subset\mathbb{R}^d$ be a
half-space. The Martin kernel with the pole at infinity for $D$ is
$s_D(x)=\delta_D^{\alpha/2}(x)$ \cite{MR2075671}. We see that
(\ref{eq:ex1_1}) and (\ref{eq:ex1_2}) hold with $C=C(d,\alpha)$ for
{\it all} $t\in(0,\infty)$ and $x,y\in{\mathbb{R}^{d}}$.
\end{exmp}

\begin{exmp}\label{ex:3}
Let $D=B^c(0,1)\subset{\mathbb{R}^{d}}$ and $d\ge \alpha$. By
the Kelvin
transform (\cite{MR2256481} or \cite{MR2365478}) and (\ref{wfg}),
% the Martin kernel is equal:
%
\begin{eqnarray*}
M_D(x)&=&
\lim_{y\rightarrow
\infty}\frac{|x|^{\alpha-d}|y|^{\alpha-d}G_B({x/|x|^2},
{y/|y|^2})}{|x_0|^{\alpha-d}|y|^{\alpha-d}G_B(
{x_0/|x_0|^2},{y/|y|^2})}
=
\frac{|x|^{\alpha-d}G_B({x/|x|^2},0)}{|x_0|^{\alpha
-d}G_B({x_0/|x_0|^2},0)},
%=\lim_{y\rightarrow
%=
\end{eqnarray*}
where
\[
G_B(z,0)={\mathcal
B}_{d,\alpha} |z|^{\alpha-d}\int^{|z|^{-2}-1}_0\frac{s^{\alpha
/2-1}}{(s+1)^{d/2}}\,ds ,\qquad
0<|z|<1.
\]
Thus, there is $c=c(x_0,d,\alpha)$ such that
%
%e39 ###
\begin{equation}\label{eq:ex2_1}
M_D(x)
%=\frac{\int^{|x|^{2}-1}_0\frac{s^{\alpha/2-1}}{(s+1)^{d/2}}ds}{
=
c
\int^{|x|^{2}-1}_0\frac{s^{\alpha/2-1}}{(s+1)^{d/2}}\,ds,\qquad
|x|\ge
1 .
\end{equation}
If $d>\alpha$, then $s_D(x)\approx1\wedge
\delta_D^{\alpha/2}(x)$, $s_D(A_{t^{1/\alpha}}(x))\approx1\wedge
(t^{1/\alpha}\vee\delta_D(x))^{\alpha/2}$, thus,
%
%e40 ###
\begin{equation}\label{eq:ex2_2}
P^x(\tau_{D}>t)\stackrel{C}{\approx}
\frac{1\wedge\delta_D^{\alpha/2}(x)}{1\wedge(t^{1/\alpha}\vee
\delta_D(x))^{\alpha/2}}=
1 \wedge\frac{\delta_D^{\alpha/2}(x)}{(1\wedge t^{1/\alpha})
^{\alpha/2}}
\end{equation}
and
\begin{eqnarray*}
\label{eq:ex2_3}
p_{D}(t, x, y)
\stackrel{C}{\approx}
\biggl(1 \wedge\frac{\delta_D^{\alpha/2}(x)}{1\wedge t^{1/2}}\biggr)
p(t,x,y)
\biggl(1 \wedge\frac{\delta_D^{\alpha/2}(y)}{1\wedge t^{1/2}}
\biggr)
\end{eqnarray*}
for all $0<t<\infty$ and $x,y\in{\mathbb{R}^{d}}$. Here
$C=C(d,\alpha)$.

For $\alpha=d=1$, (\ref{eq:ex2_1}) yields
$s_D(x)\approx\log(1+\delta_D^{1/2}(x))$,
$s_D(A_{t^{1/\alpha}}(x))\approx\log(1+(t\vee\delta_D(x))^{1/2}),$ thus,
for all $0<t<\infty$ and $x,y\in{\mathbb{R}^{d}}$ we have
%. By Theorem~\ref{theorem:bLd} with $\kappa=1/2$ and $t < \infty$:
%
%e41 ###
\begin{equation}\label{eq:ex2_4}
P^x(\tau_{D}>t)\approx
\frac{\log(1+\delta_D^{1/2}(x))}{\log(1+ (t\vee\delta_D(x))^{1/2})}=
1 \wedge\frac{\log(1+\delta_D^{1/2}(x))}{\log(1+t^{1/2})}
\end{equation}
and
\begin{eqnarray*}
\label{eq:ex2_5}
\frac{p_{D}(t, x, y)}{p(t,x,y)}
\approx
\biggl(1 \wedge\frac{\log(1+\delta_D^{1/2}(x))}{\log
(1+t^{1/2})}\biggr)
\biggl(1 \wedge\frac{\log(1+\delta^{1/2}(y))}{\log
(1+t^{1/2})}\biggr) .
\end{eqnarray*}
Sharp explicit estimates for $p_{B^c(0,R)}$ with arbitrary $R>0$
follow by scaling.
\end{exmp}

\begin{exmp}\label{ex:4}
Let $D=B^c(0,1) \subset{\mathbb{R}^{d}}$ and $1=d<\alpha$. We have that
\begin{eqnarray*}
G_{\{0\}^c}(x, y)= G_{D}(x, y)+E^xG_{\{0\}^c}(X_{T_B}, y) .
\end{eqnarray*}
Let $c_\alpha=[-2\Gamma(\alpha)\cos(\pi\alpha/2)]^{-1}$.
By \cite{MR2256481}, Lemma~4, for $x,y\in\mathbb{R}$,
\[
G_{\{0\}^c}(x,
y)=c_\alpha(|y|^{\alpha-1}+|x|^{\alpha-1}-|y-x|^{\alpha
-1}) .
\]
%
% Utilizing the integral $\int_u^v(\alpha-1)t^{\alpha-2}dt=v^{
% obtain the estimate
% \begin{equation}\label{est1}
% (1-2^{\alpha-2})c_\alpha(|x|^{\alpha-1}\wedge|y|^{\alpha-1})\le G_{
% \le2 c_\alpha(|x|^{\alpha-1}\wedge|y|^{\alpha-1})  .
% \end{equation}
If follows that
\[
G_{D}(x, y)=c_\alpha\bigl(|x|^{\alpha-1}-|x-y|^{\alpha
-1}-E^x(|X_{\tau_D}|^{\alpha-1}-|X_{\tau_D}-y|^{\alpha-1})\bigr).
\]
Since $|X_{\tau_D}|\le1$ a.s.,
$\lim_{y\rightarrow\infty}(-|x-y|^{\alpha-1}+
E^x|X_{\tau_D}-y|^{\alpha-1})=0$, for every $x\in\mathbb{R}$. If
$|x|\ge 2$, then we can find $c=c(\alpha,x_0)$ such that
\begin{eqnarray*}
M_D(x)&=&\frac{|x|^{\alpha-1}-E^x|X_{\tau
_D}|^{\alpha-1}}{|x_0|^{\alpha-1}-E^{x_0}|X_{\tau_D}|^{\alpha-1}}
=c(|x|^{\alpha-1}-E^x|X_{\tau_D}|^{\alpha-1})
\\
&\approx&
|x|^{\alpha-1}\approx\delta_D(x)^{\alpha-1}.
\end{eqnarray*}
On the other hand, by BHP, $M_D(x)\approx
\delta_D^{\alpha/2}(x)$ if $\delta(x)\le1$ (compare
Example~\ref
{tail_point}).
We thus have
$s_D(x)\approx\delta_D^{\alpha-1}(x)\wedge\delta_D^{\alpha/2}(x)$,
$s_D(A_{t^{1/\alpha}}(x))\approx
(t^{1/\alpha}\vee\delta_D(x))^{\alpha-1}\wedge
(t^{1/\alpha}\vee\delta_D(x))^{\alpha/2}$,
%. Hence by
%Theorem~\ref{theorem:bLd} with $\kappa=1/2$
and for all $0<t<\infty$, $x,y\in{\mathbb{R}^{d}}$, we obtain
%$t < \infty$:
%
%e42 ###
\begin{equation}\label{eq:ex3_1}
P^x(\tau_{D}>t)\stackrel{C}{\approx}
\frac{\delta_D^{\alpha-1}(x)\wedge
\delta_D^{\alpha/2}(x)}{(t^{1/\alpha}\vee\delta_D(x))^{\alpha
-1}\wedge
(t^{1/\alpha}\vee\delta_D(x))^{\alpha/2}}
%=
%1 \wedge\frac{\delta_D^{\alpha-1}(x)\wedge\delta_D^{
 ,
\end{equation}
hence,
\begin{eqnarray*}
\label{eq:ex3_2}
% \frac{p_{D}(t, x, y)}{p(t,x,y)}
p_{D}(t, x, y)
\stackrel{C}{\approx}
\biggl(1 \wedge\frac{\delta_D^{\alpha-1}(x)\wedge
\delta_D^{\alpha/2}(x)}{t^{1-1/\alpha}\wedge t^{1/2}}\biggr)
p(t,x,y)
\biggl(1 \wedge\frac{\delta_D^{\alpha-1}(y)\wedge
\delta_D^{\alpha/2}(y)}{t^{1-1/\alpha}\wedge
t^{1/2}}\biggr) .
\end{eqnarray*}
Here $C=C(\alpha)$. To estimate $p_{B^c(0,R)}$ with arbitrary $R>0$,
we use scaling.
% from
%(\ref{eq:ex2_3}), (\ref{eq:ex2_5}) and (\ref{eq:ex3_2}) by scaling (
\end{exmp}

\begin{defn}\label{defn:3}
We say that (open) $D$ is of class $C^{1,1}$ at scale $r>0$ if for
every $Q\in\partial D$ there exist balls
$B(x',r)\subset D$ and $B(x'',r)\subset D^c$ tangent at $Q$.
If~$D$ is $C^{1,1}$ at some
(unspecified) positive scale (hence also at smaller scales), then we
simply say $D$ is $C^{1,1}$.
\end{defn}

$C^{1,1}$ domains may be equivalently defined using local coordinates
\cite{luks-2009}.

\begin{rem}\label{rem:c11fat}
If $D$ is $C^{1,1}$ at scale $r$, then it is $(1/2,p)$-fat for all
$p\in(0,r]$.
\end{rem}

\begin{rem}\label{rem:bx}
Let $D$ be $C^{1,1}$ at scale $r$. Let $x\in D$, and let
$Q\in\partial D$ be such that $\delta_D(x)=|x-Q|$. Consider the
above balls $B(x',r)$ and $B(x'',r)$. If $\delta_D(x)<r$, then let
$B_x=B(x',r)$, otherwise $B_x=B(x,\delta_D(x))$. Thus,
$\delta_{B_x}(x)=\delta_D(x)$, and the radius of $B_x$ is $r\vee
\delta_D(x)$.
\end{rem}

\begin{exmp}\label{ex:c11}
We will verify (\ref{eq:CKS}) for $C^{1,1}$ domains $D$. For the proof we
initially assume that $D\neq{\mathbb{R}^{d}}$
is $C^{1,1}$ at scale $r=1$. Let $x\in D$. We adopt the notation of
Remark~\ref{rem:bx} and consider (the ball) $B_x$ and (the open
complement of a ball) $B^c(x'',1)$ tangent at $Q\in\partial D$.
Since $B_x\subset D \subset B^c(x'',1)$, we have
\[
P^x(\tau_{B_x}>1)\le P^x(\tau_{D}>1)\le P^x\bigl(\tau
_{B^c(x'',1)}>1\bigr) .
\]
Clearly,
$\delta_{B_x}(x)
=\delta_D(x)=|Q-x|=\delta_{B^c(x'',1)}(x)$.
By (\ref{eq:ex1_1}) and (\ref{eq:ex2_2})--(\ref{eq:ex3_1}),
\begin{eqnarray*}\label{eq:C1_1}
P^x(\tau_{D}>t)\approx
\biggl(1 \wedge\frac{\delta_D(x)}{t^{1/\alpha}}\biggr)^{\alpha/2}
 ,\qquad t\le1 .
\end{eqnarray*}
By Remark~\ref{rem:c11fat} and Theorem~\ref{theorem:oppz}, there is
$C=C(d,\alpha)$ such that, for all $x,y\in{\mathbb{R}^{d}}$,
\begin{eqnarray*}
\label{eq:gc11cs}
p_{D}(t, x, y)\stackrel{C}{\approx}\biggl(
1\wedge\frac{\delta_D(x)^{\alpha/2}}{t^{1/2}}
\biggr)
\biggl(\frac{t}{|x|^{d+\alpha}} \land
t^{-d/\alpha}\biggr)
\biggl(1\wedge
\frac{\delta_D(y)^{\alpha/2}}{t^{1/2}}
\biggr)
 ,\qquad t\le1 .
\end{eqnarray*}
If $D$ is $C^{1,1}$ at a scale $r<1$, then $r^{-1}D$ is
$C^{1,1}$ at scale $1$. This yields (\ref{eq:CKS}) in time range
$0<t\le r^\alpha$. Remark~\ref{comp} allows for an extension to all
$t\in(0,1]$, with a constant depending on $d$, $\alpha$ and
$r$. The case of $D={\mathbb{R}^{d}}$ is trivial.
\end{exmp}

Further estimates for $C^{1,1}$ domains will be given in
Proposition~\ref{cor:espc11}, Theorem~\ref{cor:c2} and Corollary~\ref
{cor:c4}.
% The remainder of this section is devoted to the $C^{1,1}$
% condition. We pay considerable attention to the dependance of the
% comparability constants of our estimates on geometric characteristics
% of the domains in question.

\begin{exmp}\label{ex:sL}
%The case of general $(\kappa,r)$-fat domain is obtained by scaling
%(see also Remark~\ref{rem:inc}).
%Lipschitz and $C^{1,1}$ domains are a special case.
Let $d\ge 2$. For $x=(x_1,\ldots,x_{d-1},x_d)\in{\mathbb{R}^{d}}$
we denote $\tilde{x}=(x_1,\ldots,x_{d-1})$, so that
$x=(\tilde{x},x_d)$.
Let $\lambda<\infty$. We consider a Lipschitz function
$\gamma\dvtx \mathbb{R}^{d-1}\to\mathbb{R}$, that is,
$|\gamma(\tilde{x})-\gamma(\tilde{y})|\le\lambda|\tilde
{x}-\tilde{y}|$.
We define a special Lipschitz domain
$D=\{x=(\tilde{x},x_d)\in
{\mathbb{R}^{d}}\dvtx x_{d}>\gamma(\tilde{x})\}$. For such $D$ the
geometric notions of
Theorem~\ref{theorem:bLd} become more explicit as we will see below.
We note that $D$ is $((2\sqrt{1+\lambda^2})^{-1},r)$-fat for all
$r>0$ (\cite{MR2013738}, Remark~1).
For $x=(\tilde{x},x_d)\in D$ and $r>0$ we define
$x^{(r)}=(\tilde{x},\gamma(\tilde{x})+r)$. If $x$ is close to
$\partial D$,
then $x^{(1)}$ dominates $x$ in the direction of the last coordinate.
We note that $P^{x^{(1)}}(\tau_D>1)\ge  c>0$. Here $c=c(d,\alpha
,\lambda)$.
%We fix an arbitrary reference point $x_0\in D$.
%and consider the Martin
%kernel with the pole at infinity \cite{MR2365478},
% M(x) =
% \lim_{D \ni v \rightarrow\infty} \frac{G_{D}(x, v)}{G_{D}(x_0, v)}
%   , \quad x \in\Rd .
%The limit exists, and $M$ is $\alpha$-harmonic on $D$, see
By Remark~\ref{comp} and BHP,
%
%e43 ###
\begin{equation}
\label{eq:spld}
P^{x}(\tau_D>1)\stackrel{C}{\approx} 1\wedge\frac
{M_D(x)}{M_D(x^{(1)})} ,\qquad
x\in D ,
\end{equation}
where $C=C(\alpha,d,\lambda)$.
By scaling,
% of the Green function,
the Martin kernel with the pole at infinity for $rD$ is a constant multiple
of $M_D(x/r)$. By
%this and
(\ref{eq:spld}), we obtain
%
%e44 ###
\begin{equation}\label{eq:spc11}
P^x(\tau_{D}>t)=
P^{t^{-1/\alpha}x}(\tau_{t^{-1/\alpha}D}>1)\stackrel{C}{\approx}
1\wedge\frac{M_D(x)}{M_D(x^{(t^{1/\alpha})})}
 ,\qquad x\in D .
\end{equation}
We note in passing that (\ref{eq:spc11}) agrees with (\ref{eq:spkf})
because $r\mapsto M_D(x^{(r)})$ is increasing~\cite{2002-KBTKANijm}.
Or, in our previous notation we can take $A_r(x,\kappa,D)=x^{(r\vee
(x_d-\gamma(\tilde{x}))}$.
We substitute (\ref{eq:spc11}) into (\ref{eq:get}) so that
for all $0<t<\infty$ and $x,y\in D$ (in fact, by regularity, for
$x,y\in{\mathbb{R}^{d}}$) we have
\begin{eqnarray*}
\label{eq:getsld}
p_{D}(t, x, y)
\stackrel{C}{\approx}
\biggl(1\wedge\frac{M_D(x)}{M_D(x^{(t^{1/\alpha})})}\biggr)
p(t,x,y)
\biggl(1\wedge\frac{M_D(y)}{M_D(y^{(t^{1/\alpha})})}\biggr)
 .
\end{eqnarray*}
%
%By BHP (\cite{MR2365478}), the Martin kernel may be replaced here by
%some other
%(super) harmonic functions.
% For $x\in\Rd\setminus\{0\}$ we denote by $\theta(x)$ the angle
% between $x$ and the point $(0,\ldots,0,1)\in\Rd$.
% We fix $0<\Theta<\pi$ and consider the right circular cone
% $V=\{x\in\Rd\setminus\{0\}:  \theta(x)<\Theta\}$.
% Clearly, $rV=V$ for every $r>0$.
\end{exmp}

\begin{exmp}\label{ex:6}
For {circular cones}
$V$ \cite{2009KBTGcm}
%by Example~\ref{ex:sL}
we have
%
%e45 ###
\begin{equation}
\label{eq:sMk}
M_V(x)=|x|^\beta M_V(x/|x|) ,\qquad x\neq0 ,
\end{equation}
where $0\le\beta<\alpha$ is a characteristic of the cone; see
\cite{MR2075671}.
% , \cite{MR2182071}).
By \cite{MR2213639}, Lemma~3.3,
\begin{eqnarray*}
\label{eq:ojmcc}
M_V(x)\approx\delta_V(x)^{\alpha/2}|x|^{\beta-\alpha/2} ,\qquad
x\in{\mathbb{R}^{d}} ;
\end{eqnarray*}
see also \cite{2009KBTGcm} and \cite{MR2213639}.
%Here $\delta(x)=\delta_{D}(x)$.
Considering (\ref{eq:spld}), by simple manipulations, we obtain
%
%e46 ###
\begin{equation}\label{eq:sscone}
1\wedge\frac{\delta_V(x)^{\alpha/2}|x|^{\beta-\alpha/2}
}
{(1 \vee|\tilde{x}|)^{\beta-\alpha/2}}
\stackrel{C}{\approx}
\bigl(1 \wedge\delta_V(x)^{\alpha/2}\bigr)(1 \wedge|x|)^{\beta-\alpha/2} ,
\end{equation}
where $C=C(\lambda)$. By (\ref{eq:get}) and scaling, we get (\ref{eq:etdG2}).
\end{exmp}

The interested reader may find more references on stable processes and
Brownian motion in
cones in \cite{2009KBTGcm}. Note that (\ref{eq:sMk}) holds for
generalized open cones, that is, open sets $\varnothing\neq V\subset
{\mathbb{R}^{d}}$
such that $kV=V$ for all $k>0$ \cite{MR2075671}.

\begin{exmp}\label{ex:8}
Let $d=1,2,\ldots$ and $V={\mathbb{R}^{d}}\setminus\{x_d=0\}$.
This generalized cone is {non-Lipschitz} but it is $(1/2,r)$-fat for
every $r>0$.
Let $1<\alpha<2$.
From \cite{MR2075671}, Example 3.3, we have $M_V(x)=|x_d|^{\alpha-1}$
(the decay near a hyperplane is slower than near a half-space).
We consider $t=1$ in (\ref{eq:ppkf}).
We let $A_1(x)=(\tilde{x}, x_d+1/2)$ if $x_d>0$ and $A_1(x)=(\tilde
{x},x_d-1/2)$ otherwise. Thus,
\begin{eqnarray*}
\frac{M_V(x)}{M_V(A_1(x))} =
\frac{|x_d|^{\alpha-1}}{(|x_d|+1/2)^{\alpha-1}}
\approx(1 \wedge|x_d|)^{\alpha-1} .
\end{eqnarray*}
%
%This is similar to (\ref{eq:sscone}).
By (\ref{eq:get}) and scaling, we obtain the following analogue of
(\ref{eq:etdG2}):
%
%e47 ###
\begin{equation}\label{eq:etdGgc}
\quad\ \frac{p_{V}(t,x,y)}{p(t,x,y)}
\approx
\biggl(1\wedge\frac{\delta_V(x)}{t^{1/\alpha}}\biggr)^{\alpha-1}
\biggl(1\wedge\frac{\delta_V(y)}{t^{1/\alpha}}\biggr)^{\alpha-1}
, \qquad t>0 ,x,y\in{\mathbb{R}^{d}} .
\end{equation}
We note that $V$ is the complement of a point if $d=1$.
%It is also possible to give an explicit estimate for the complement of
%a half-line on %the plane, and its higher dimensional counterparts.
%The interested reader is referred to \cite[Theorem 3.2
%(iii)]{MR2213639}.
\end{exmp}

If $D$ is bounded and $\kappa>0$ is fixed, then $D$ is {\it not}\/
$(\kappa,r)$ at large scales
$r$, and the asymptotics of the probability of survival are
exponential. Indeed, for the fractional Laplacian with Dirichlet
condition on $D^c$
we let $\lambda_1>0$ be its first eigenvalue and
$\phi_1>0$ the corresponding eigenfunction [normalized in $L^2(D,dx)$];
see \cite{MR1643611}.
The following approximation results from the intrinsic
ultracontractivity of
% the killed semigroup of the fractional Laplacian
%in
every bounded domain \cite{MR1643611}:
\begin{eqnarray*}
\label{eq:ppkflt}
p_{D}(t, x, y)
\approx\phi_1(x)\phi_1(y)e^{-\lambda_1 t} , \qquad
t\ge 1 ,  x,y\in{\mathbb{R}^{d}} .
\end{eqnarray*}
Here comparability constants depend on $D$ and $\alpha$ (see also
Proposition~\ref{cor:espc11} below).
Given that infinity is {\it inaccessible} \cite{MR2365478} from
bounded $D$, it is of considerable interest to understand
the behavior of the heat kernel related to accessible and
inaccessible points of $D$ (see also \cite{2009-MK-pa} in this connection).

In the remainder of the paper we will study
$C^{1,1}$ domains in more detail.
We focus on unbounded domains, large times and dependence of the
comparability constants on global geometry of the domains.
%We denote $\delta_U(x)=\operatorname{dist}(x,U^c)$, the distance of $x$ to the
%complement of $U$.
%For instance we have
% There is $C_2=C_2(\alpha,d)$ such that for all $r\ge 1$ and $x\in
% \Rd$,
% \begin{equation}\label{ball_tail}
% P^x(\tau_{B(0,r)}>1)\appc{C_2}
% \delta_{B(0,r)}^{\alpha/2}(x)\wedge1 ,
% \end{equation}\
% and
% \begin{equation}\label{outer-ball-tail}
% P^x(T_{B(0,r)}>1) \appc{C_2}
% \delta_{B(0,r)^c}^{\alpha/2}(x)\wedge1 .
% \end{equation}
% Let $x\in B=B(0,r)$. If $\delta_{B}(x)\ge 1/4$ then $1\ge
% P^x(\tau_{B}>1)\ge  P^0(\tau_{B(0,1/4)}>1)>0$, which yields
% (\ref{ball_tail}). If $\delta_B(x)<1/4$ then we let
% $A=x[1-1/(2|x|)]$ and $U= B\cap B(x,|x-A|+1/6)$,
% according to Definition~\ref{def:U} and Remark~\ref{rem:ball} (we
% consider $D=B$ in Definition~\ref{def:U}).
% The function $y\mapsto P^y(X_{\tau_{U}}\in B)$ is harmonic in $U$,
% and $P^A(X_{\tau_U}\in B)\ge  c>0$. Inspecting (\ref{wfg}) we see
% that $r^{\alpha/2}G_{B}(y,0)\approx\delta_B(y)^{\alpha/2}$ if
% $\delta_B(y)<1/4$. By BHP, $P^x(X_{\tau_{U}}\in B)\approx
% \delta_B(x)^{\alpha/2}$, too. This proves (\ref{ball_tail}), see
% Remark~\ref{comp}. The proof of (\ref{outer-ball-tail}) is similar
% (see \cite{MR0126885} for the Green function of the complement of the
% ball).

Example~\ref{ex:1} and intrinsic ultracontractivity yield the following result.
\begin{lem}\label{hittingtime}
There exist $\lambda_1 =\lambda_1 (\alpha, d)>0$ and
$C=C(\alpha,d)$ such that for all $r>0$, $t>0$ and $x\in{\mathbb
{R}^{d}}$ we have
\begin{eqnarray*}
\label{eq:*}
P^x\bigl(\tau_{B(0,r)}>t\bigr)\stackrel{C}{\approx} \biggl[1 \wedge
\biggl(\frac{\delta
_{B(0,r)}(x)}{r\wedge t^{1/\alpha}}\biggr)^{\alpha/2}
\biggr]e^{-\lambda_1 t/r^\alpha} .
\end{eqnarray*}
%
% There are constants $C_8=C_8(\alpha,d)$ and $C_9=C_9(\alpha,d)$ such
%that
% \begin{equation}
% \label{eq:**}
% C_8 [(\frac{\delta_{B^c(0,r)}(x)}{ t^{1/\alpha}})^{
% [(\frac{\delta_{B^c(0,r)}(x)}{r\wedge t^{1/\alpha}}
%)^{\alpha/2}\wedge1] .
% \end{equation}
% Moreover for $d>\alpha$ we have
% \begin{equation}
% \label{eq:***}
% P^x(T_{B(0,r)}>t)\appc{C_{10}}
% (\frac{\delta_{B^c(0,r)}(x)}{r\wedge t^{1/\alpha}})^{
% \end{equation}
% where $C_{10}=C_{10}(\alpha,d)$.
\end{lem}

\begin{lem}\label{hittingtime1}
Let $d>\alpha$, $0<r<R$, $W=B(0,r)\cup B^c(0,R)$. There is
$c=c(\alpha,d)$ such that for all $t>0$ and $x\in{\mathbb{R}^{d}}$
we have
\begin{eqnarray*}
\label{eq:****}
P^x(\tau_{W}>t)\ge  c \biggl(\frac rR\biggr)^\alpha
\biggl[1\wedge\biggl(\frac{\delta_{B(0,r)}(x)}{r\wedge
t^{1/\alpha}}\biggr)^{\alpha/2}\biggr] .
\end{eqnarray*}
\end{lem}

\begin{pf}
By scaling, we only need to consider $r=1< R$. By \cite{MR0126885}, we obtain
\begin{eqnarray*}
% \label{eq:uu}
P^x\bigl(T_{B(0,1)}=\infty\bigr)&=&
\frac{\Gamma(d/2)}{\Gamma((d-\alpha)/2)\Gamma(\alpha/2)}
\int_0^{|x|^2-1}
\frac{u^{\alpha/2-1}}{(u+1)^{d/2}}
\,du
\\
&\approx&
1 \wedge\delta_{B^c(0,1)}^{\alpha/2}(x)
\end{eqnarray*}
%
%we
% have
%$$ P^y(T_{B(0,R)}=\infty)\approx
%(\frac{\delta_{B^c(0,R)}(y)}R)^{\alpha/2}\wedge1 ,\quad
%y\in\Rd ,$$ thus
[compare (\ref{eq:ex2_1})]. Thus, there is $c=c(d,\alpha)$ such that
\[
P^y\bigl(T_{B(0,R)}=\infty\bigr)\ge  c>0,\qquad |y|>2R.
\]
Let $x\in B(0,1)$. For $t\ge 1$ we use (\ref{eq:tpk}) to obtain
\begin{eqnarray*}
P^x(\tau_{W}>t)&\ge &
P^x(\tau_{W}=\infty)
\\
&\ge & E^x\bigl\{\bigl|X_{\tau_{B(0,1)}} \bigr|\ge 2R; P^{X_{\tau
_{B(0,1)}}
}\bigl(T_{B(0, R)}=\infty\bigr)\bigr\}
\\
&\ge & cP^x\bigl(\bigl|X_{\tau_{B(0,1)} }\bigr|\ge 2R\bigr)
\ge  c \frac{1}{R^\alpha}\delta^{\alpha/2}_{B(0,1)}(x) .
\end{eqnarray*}
By (\ref{eq:ex1_1}), for $t\le1$
we even have
\begin{eqnarray*}
P^x(\tau_{W}>t)&\ge & P^x\bigl(\tau_{B(0,1)}>t\bigr) \approx
1\wedge\biggl(\frac{\delta_{B(0,1)}(x)}{1\wedge t^{1/\alpha
}}\biggr)^{\alpha/2}.
\end{eqnarray*}
\upqed\end{pf}

The $C^{1,1}$ condition at a given scale fails to determine the
{fatness} of $D$
at larger scales and, consequently, the exact asymptotics of the
survival probability.
The following is a substitute.
% and the first (general) estimate in
%Proposition~\ref{cor:espc11} cannot be sharp. The (second) estimate
%for $C^{1,1}$
%domains of bounded complement is sharp except for the factor
%$r/R$, which controls global geometry of $D$.

\begin{prop}\label{cor:espc11}
If $D$ is $C^{1,1}$ at some scale $r>0$, then
\begin{eqnarray}
\label{eq:gbc11}
&&C^{-1} e^{- {\lambda_1 t/(r\vee\delta_D(x))^\alpha}}
\biggl[1\wedge\biggl(\frac{\delta_D(x)}{r\wedge t^{1/\alpha
}}\biggr)^{\alpha/2}\biggr]\nonumber
\\[-8pt]\\[-8pt]
&&\qquad \le P^x(\tau_{D}>t)
\le C\biggl[1\wedge\biggl(\frac{\delta_D(x)}{r\wedge
t^{1/\alpha
}}\biggr)^{\alpha/2}\biggr]\nonumber
\end{eqnarray}
for all $t>0$ and $x\in{\mathbb{R}^{d}}$.
Here $C=C(\alpha,d)$
and $\lambda_1 =\lambda_1 (\alpha,d)$.

If {\rm also} $d>\alpha$ and ${\rm
diam}(D^c)<\infty$, then for all $t>0$ and $x\in{\mathbb{R}^{d}}$,
%
%e48 ###
\begin{equation}
\label{eq:sbc11}
P^x(\tau_{D}>t)\ge
C^{-1} \biggl(\frac{r}{\operatorname{diam}(D^c)}\biggr)^{\alpha}
\biggl[1\wedge\biggl(\frac{\delta_D(x)}{r\wedge
t^{1/\alpha}}\biggr)^{\alpha/2}\biggr] .
\end{equation}
\end{prop}

\begin{pf}
Consider $x\in D$, $B_x\subset D$
and $B(x'',r)\subset D^c$ of Remark~\ref{rem:bx}.
Clearly, $\tau_{B_x}\le\tau_{D}\le T_{B(x'',r)}$, thus,
\[
P^x(\tau_{B_x}>t)\le P^x(\tau_{D}>t)\le
P^x\bigl(T_{B(x'',r)}>t\bigr) .
\]
Lemma \ref{hittingtime} yields the estimate
\begin{eqnarray*}C^{-1} e^{-\lambda_1
t/(r\vee\delta_D(x)) ^\alpha} \biggl[1\wedge\biggl(\frac
{\delta_D(x)}{(r\vee\delta_D(x))\wedge t^{1/\alpha}}
\biggr)^{\alpha/2}\biggr]&\le&P^x(\tau_{D}>t)
\end{eqnarray*}
and
\[
P^x(\tau_{D}>t)\le C \biggl[1\wedge\biggl(\frac{\delta
_D(x)}{r\wedge
t^{1/\alpha}}\biggr)^{\alpha/2}\biggr] ,
\]
which simplifies to (\ref{eq:gbc11}) as
$\delta_D(x)>r$ yields
$\delta_D(x)/[(r\vee\delta_D(x))\wedge t^{1/\alpha}]\ge 1$.
To prove (\ref{eq:sbc11}), we consider $\rho=\operatorname{diam}(D^c)\ge 2r$,
the center, say, $x_0$, of $B_x$, and $W:=B_x\cup B^c(x_0,\rho+r\vee
\delta_D(x))\subset D$.
By Lemma \ref{hittingtime1} and Remark~\ref{rem:bx},
\begin{eqnarray*} P^x(\tau_{D}>t)&\ge & P^x(\tau_{W}>t)
\\
&\ge &
c
\biggl(\frac{r\vee
\delta_D(x)}{\rho+r\vee\delta_D(x)}\biggr)^\alpha\biggl[1\wedge
\biggl(\frac
{\delta_D(x)}{(r\vee\delta_D(x))\wedge
t^{1/\alpha}}\biggr)^{\alpha/2}\biggr]
\\
&\ge & c
\biggl(\frac r\rho\biggr)^{\alpha} \biggl[1\wedge\biggl(\frac
{\delta_D(x)}{r\wedge
t^{1/\alpha}}\biggr)^{\alpha/2}\biggr] .
\end{eqnarray*}
\upqed\end{pf}

In view of Theorem~\ref{theorem:oppz}, (\ref{eq:gbc11}) mildly strengthens
\cite{CKS2008}, Theorem~1.1(i) [i.e., (\ref{eq:CKS}) above]. We also
get the following result.

%The proof of the lower bound of the estimate for bounded $D^c$ is
%completed.
%If $D$ is $C^{1,1}$ at scale $r>0$ then
%Theorem~\ref{theorem:ext} and Proposition~\ref{cor:espc11} yield
%$C=C(d,\alpha,r)$ such that for
%all $x,y\in D$ we have
% \label{eq:gc11cs}
%p_{D}(t, x, y)\appc{C}(
%(
%t^{-d/\alpha}) ,\quad t\le1 .
% \end{equation}
%We thus reproved
%For $D=H$, the halfspace, (\ref{eq:gc11cs}) holds for all $t>0$ and
%$x,y\in\Rd$ with
%$C=C(d,\alpha)$, see Example~\ref{tail_point}.
%In what follows we will study arbitrary $t>0$.

\begin{theorem}\label{cor:c2}
Let $d>\alpha$. If $D$ is $C^{1,1}$ at scale $r$ and ${\rm
diam}(D^c)<\infty$, then
\begin{eqnarray*}
\label{eq:gc11}
&&C^{-1}\biggl(\frac{r}{\operatorname{diam}(D^c)}\biggr)^{2\alpha}
\\
&&\qquad \le
\frac{p_{D}(t, x, y)}{[1\wedge({\delta
_D(x)/(r\wedge
t^{1/\alpha})})^{\alpha/2}]
p(t,x,y)
[1\wedge
({\delta_D(y)/(r\wedge
t^{1/\alpha})})^{\alpha/2}
]}
\\
&&{}\qquad \le C
\end{eqnarray*}
for all $t>0$ and $x,y\in{\mathbb{R}^{d}}$. Here $C=C(\alpha,d)$.
%, and, in particular, for $x,y\in B^c$ we have
% \label{eq:gb}
% p_{B^c(0,1)}(t, x, y)
%[ (\frac
% {|x|-1}{1\wedge
% t^{1/\alpha}})^{\alpha/2}\wedge1][
% (\frac{|y|-1}{1\wedge
% t^{1/\alpha}})^{\alpha/2}\wedge1] p(t, x,
% y) .
\end{theorem}

\begin{pf}
The result follows from (\ref{eq:sbc11}) and Corollary~\ref{theorem:ext}.
\end{pf}

A similar result (with less control of the constants) is given in
\cite{2009ZCJT}.\footnote{Paper \cite{2009ZCJT} appeared on arXiv after
the first draft \cite{2009KBTGMR} of the present paper.}

\begin{rem} We consider the recurrent case $\alpha\ge  d=1$.
If $D\subset\mathbb{R}$ is the complement of a finite union of
bounded closed intervals, then
\begin{eqnarray*}
P^x(\tau_D>t)&\stackrel{C}{\approx}&1\wedge
\frac{\delta_D(x)^{\alpha-1}\wedge
\delta_D(x)^{\alpha/2}}{t^{1-1/\alpha}\wedge
t^{1/2}} ,
\qquad t>0 , x\in{\mathbb{R}^{d}} ,
\mbox{ if  $\alpha>1$,}
\\
P^x(\tau_D>t)&\stackrel{C}{\approx}&1\wedge\frac{\log(1+ \delta
_D(x)^{1/2})}{\log(1+
t^{1/2})} ,
\qquad t>0 , x\in{\mathbb{R}^{d}} ,
\mbox{ if  $\alpha=1$,}
\end{eqnarray*}
where $C=C(D,\alpha)$.
%We omit t
The estimates
%since it
follow easily from Examples \ref{tail_point} and \ref{ex:3}.
\end{rem}

\begin{cor}\label{cor:c4}
If $D\subset\mathbb{R}$ is the complement to a finite union of
bounded closed intervals, then
$C=C(D,\alpha)$ exists such that for all $t>0$ and $x,y\in\mathbb{R}$,
\begin{eqnarray*} \frac{p_{D}(t, x, y)}{p(t, x,
y)}\stackrel{C}{\approx}
\biggl[1\wedge\frac{\delta_D(x)^{\alpha-1}\wedge
\delta_D(x)^{\alpha/2}}{t^{1-1/\alpha}\wedge
t^{1/2}}\biggr]\biggl[1\wedge
\frac{\delta_D(y)^{\alpha-1}\wedge
\delta_D(y)^{\alpha/2}}{t^{1-1/\alpha}\wedge
t^{1/2}}\biggr]\vspace*{-1pt}
\end{eqnarray*}
for $\alpha>1$, while for $\alpha=1$ we have
\begin{eqnarray*} \frac{p_{D}(t, x, y)}{p(t, x,
y)}\stackrel{C}{\approx} \biggl[1\wedge\frac{\log(1+
\delta_D(x)^{1/2})}{\log(1+ t^{1/2})}\biggr]\biggl[1\wedge
\frac{\log(1+\delta_D(y)^{1/2})}{\log(1+
t^{1/2})}\biggr].\vspace*{-2pt}
\end{eqnarray*}

% Let $d=1$. For $D=[-R,R]^c$ and $\alpha>1$ we have

% y)}&\approx&
% [\frac{(|x|-R)^{\alpha-1}\wedge
% R^{\alpha/2-1}(|x|-R)^{\alpha/2}}{t^{1-1/\alpha}\wedge
% R^{\alpha/2-1} t^{1/2}}\wedge1][
% \frac{(|y|-R)^{\alpha-1}\wedge
% R^{\alpha/2-1}(|y|-R)^{\alpha/2}}{t^{1-1/\alpha}\wedge
% R^{\alpha/2-1} t^{1/2}}\wedge1] .\end{eqnarray*}

% For $\alpha=1$ we have

% y)}&\approx& [\frac{\log(1+
% (|x-R|/R)^{1/2})}{\log(1+ (t/R)^{1/2})}\wedge1][
% \frac{\log(1+ (|x-R|/R)^{1/2})}{\log(1+
% (t/R)^{1/2})}\wedge1].\end{eqnarray*}
\end{cor}

\section*{Acknowledgments} Results of the paper were reported
at the workshop {\it Schr\"odinger operators and stochastic processes},
Wroc\l{}aw, 14--15 May 2009, and at
{\it The Third International Conference on Stochastic Analysis and Its
Applications}, Beijing, 13--17 July 2009. The authors are grateful to
the organizers
for this opportunity. We also thank an anonymous referee for a question
which led to our conjecture in the \hyperref[sec1]{Introduction} about (\ref{eq:get})
for irregular boundary points of $D$.

%dla kompilacji w bibtexu uaktualnionej bazy bibliograficznej
%odkomentuj kolejne dwie linie oraz skasuj ponizsza bibliografie az do
%...a po kompilacji w Latexu i Bibtexu zakomentu powyzsze dwie linie,
%skopiuj ponizej zawartosc zbioru.bbl i dopisz \end{document}

 \def\cprime{$'$}

\printaddresses


\begin{thebibliography}{10}

%b1 ###
\bibitem{MR2075671}
\begin{barticle}[mr]
\bauthor{\bsnm{Ba{\~n}uelos},~\bfnm{Rodrigo}\binits{R.}} \AND
  \bauthor{\bsnm{Bogdan},~\bfnm{Krzysztof}\binits{K.}}
(\byear{2004}).
\btitle{Symmetric stable processes in cones}.
\bjournal{Potential Anal.}
\bvolume{21}
\bpages{263--288}.
\bid{doi={10.1023/B:POTA.0000033333.72236.dc}, mr={2075671}}
\end{barticle}
\endbibitem

%b2 ###
\bibitem{MR2438694}
\begin{barticle}[mr]
\bauthor{\bsnm{Ba{\~n}uelos},~\bfnm{Rodrigo}\binits{R.}} \AND
  \bauthor{\bsnm{Kulczycki},~\bfnm{Tadeusz}\binits{T.}}
(\byear{2008}).
\btitle{Trace estimates for stable processes}.
\bjournal{Probab. Theory Related Fields}
\bvolume{142}
\bpages{313--338}.
\bid{doi={10.1007/s00440-007-0106-x}, mr={2438694}}
\end{barticle}
\endbibitem

%b3 ###
\bibitem{MR2492992}
\begin{barticle}[mr]
\bauthor{\bsnm{Barlow},~\bfnm{Martin~T.}\binits{M.~T.}},
  \bauthor{\bsnm{Grigor'yan},~\bfnm{Alexander}\binits{A.}} \AND
  \bauthor{\bsnm{Kumagai},~\bfnm{Takashi}\binits{T.}}
(\byear{2009}).
\btitle{Heat kernel upper bounds for jump processes and the first exit time}.
\bjournal{J. Reine Angew. Math.}
\bvolume{626}
\bpages{135--157}.
\bid{doi={10.1515/CRELLE.2009.005}, mr={2492992}}
\end{barticle}
\endbibitem

%b4 ###
\bibitem{MR0119247}
\begin{barticle}[mr]
\bauthor{\bsnm{Blumenthal},~\bfnm{R.~M.}\binits{R.~M.}} \AND
  \bauthor{\bsnm{Getoor},~\bfnm{R.~K.}\binits{R.~K.}}
(\byear{1960}).
\btitle{Some theorems on stable processes}.
\bjournal{Trans. Amer. Math. Soc.}
\bvolume{95}
\bpages{263--273}.
\bid{mr={0119247}}
\end{barticle}
\endbibitem

%b5 ###
\bibitem{MR0126885}
\begin{barticle}[mr]
\bauthor{\bsnm{Blumenthal},~\bfnm{R.~M.}\binits{R.~M.}},
  \bauthor{\bsnm{Getoor},~\bfnm{R.~K.}\binits{R.~K.}} \AND
  \bauthor{\bsnm{Ray},~\bfnm{D.~B.}\binits{D.~B.}}
(\byear{1961}).
\btitle{On the distribution of first hits for the symmetric stable processes.}
\bjournal{Trans. Amer. Math. Soc.}
\bvolume{99}
\bpages{540--554}.
\bid{mr={0126885}}
\end{barticle}
\endbibitem

%b6 ###
\bibitem{MR1741527}
\begin{barticle}[mr]
\bauthor{\bsnm{Bogdan},~\bfnm{Krzysztof}\binits{K.}}
(\byear{2000}).
\btitle{Sharp estimates for the {G}reen function in {L}ipschitz domains}.
\bjournal{J. Math. Anal. Appl.}
\bvolume{243}
\bpages{326--337}.
\bid{doi={10.1006/jmaa.1999.6673}, mr={1741527}}
\end{barticle}
\endbibitem

%b7 ###
\bibitem{MR1671973}
\begin{barticle}[mr]
\bauthor{\bsnm{Bogdan},~\bfnm{Krzysztof}\binits{K.}} \AND
  \bauthor{\bsnm{Byczkowski},~\bfnm{Tomasz}\binits{T.}}
(\byear{1999}).
\btitle{Potential theory for the {$\alpha$}-stable {S}chr\"odinger operator on
  bounded {L}ipschitz domains}.
\bjournal{Studia Math.}
\bvolume{133}
\bpages{53--92}.
\bid{mr={1671973}}
\end{barticle}
\endbibitem

%b8 ###
\bibitem{MR1825645}
\begin{barticle}[mr]
\bauthor{\bsnm{Bogdan},~\bfnm{Krzysztof}\binits{K.}} \AND
  \bauthor{\bsnm{Byczkowski},~\bfnm{Tomasz}\binits{T.}}
(\byear{2000}).
\btitle{Potential theory of {S}chr\"odinger operator based on fractional
  {L}aplacian}.
\bjournal{Probab. Math. Statist.}
\bvolume{20}
\bpages{293--335}.
\bid{mr={1825645}}
\end{barticle}
\endbibitem

%b9 ###
\bibitem{KB-TB-MR-TK-RS-ZV}
\begin{bbook}[mr]
\bauthor{\bsnm{Bogdan},~\bfnm{Krzysztof}\binits{K.}},
  \bauthor{\bsnm{Byczkowski},~\bfnm{Tomasz}\binits{T.}},
  \bauthor{\bsnm{Kulczycki},~\bfnm{Tadeusz}\binits{T.}},
  \bauthor{\bsnm{Ryznar},~\bfnm{Michal}\binits{M.}},
  \bauthor{\bsnm{Song},~\bfnm{Renming}\binits{R.}} \AND
  \bauthor{\bsnm{Vondra{\v{c}}ek},~\bfnm{Zoran}\binits{Z.}}
(\byear{2009}).
\btitle{Potential Analysis of Stable Processes and Its Extensions}
(\beditor{P.~Graczyk and A. Stos}, eds.).
\bseries{Lecture Notes in Math.}
\bvolume{1980}.
\bpublisher{Springer}, \baddress{Berlin}.
\bid{mr={2569321}}
\end{bbook}
\endbibitem

%b10 ###
\bibitem{2009KBTGcm}
\begin{barticle}[auto:SpringerTagBib|2009-01-14|16:51:27]
\bauthor{\bsnm{Bogdan},~\bfnm{K.}\binits{K.}} \AND
  \bauthor{\bsnm{Grzywny},~\bfnm{T.}\binits{T.}}
  (\byear{2010}).
\btitle{Heat kernel of fractional {L}aplacian in cones}.
  \bjournal{Colloq. Math.}
\bvolume{118}
\bpages{365--377}.
\end{barticle}
\endbibitem

%b11 ###
\bibitem{2009KBTGMR}
\begin{bmisc}[auto:SpringerTagBib|2009-01-14|16:51:27]
\bauthor{\bsnm{Bogdan},~\bfnm{K.}\binits{K.}},
  \bauthor{\bsnm{Grzywny},~\bfnm{T.}\binits{T.}} \AND
  \bauthor{\bsnm{Ryznar},~\bfnm{M.}\binits{M.}}
  (\byear{2009}).
\btitle{Heat kernel estimates for the fractional {L}aplacian}.
\bnote{Preprint. Available at} \url{http://arxiv.org/abs/0905.2626v1}.
\end{bmisc}
\endbibitem

%b12 ###
\bibitem{2008KBTJ}
\begin{barticle}[mr]
\bauthor{\bsnm{Bogdan},~\bfnm{Krzysztof}\binits{K.}},
  \bauthor{\bsnm{Hansen},~\bfnm{Wolfhard}\binits{W.}} \AND
  \bauthor{\bsnm{Jakubowski},~\bfnm{Tomasz}\binits{T.}}
(\byear{2008}).
\btitle{Time-dependent {S}chr\"odinger perturbations of transition densities}.
\bjournal{Studia Math.}
\bvolume{189}
\bpages{235--254}.
\bid{doi={10.4064/sm189-3-3}, mr={2457489}}
\end{barticle}
\endbibitem

%b13 ###
\bibitem{MR2283957}
\begin{barticle}[mr]
\bauthor{\bsnm{Bogdan},~\bfnm{Krzysztof}\binits{K.}} \AND
  \bauthor{\bsnm{Jakubowski},~\bfnm{Tomasz}\binits{T.}}
(\byear{2007}).
\btitle{Estimates of heat kernel of fractional {L}aplacian perturbed by
  gradient operators}.
\bjournal{Comm. Math. Phys.}
\bvolume{271}
\bpages{179--198}.
\bid{doi={10.1007/s00220-006-0178-y}, mr={2283957}}
\end{barticle}
\endbibitem

%b14 ###
\bibitem{MR2365478}
\begin{barticle}[mr]
\bauthor{\bsnm{Bogdan},~\bfnm{Krzysztof}\binits{K.}},
  \bauthor{\bsnm{Kulczycki},~\bfnm{Tadeusz}\binits{T.}} \AND
  \bauthor{\bsnm{Kwa{\'s}nicki},~\bfnm{Mateusz}\binits{M.}}
(\byear{2008}).
\btitle{Estimates and structure of {$\alpha$}-harmonic functions}.
\bjournal{Probab. Theory Related Fields}
\bvolume{140}
\bpages{345--381}.
\bid{doi={10.1007/s00440-007-0067-0}, mr={2365478}}
\end{barticle}
\endbibitem

%b15 ###
\bibitem{2002-KBTKANijm}
\begin{barticle}[mr]
\bauthor{\bsnm{Bogdan},~\bfnm{K.}\binits{K.}},
  \bauthor{\bsnm{Kulczycki},~\bfnm{T.}\binits{T.}} \AND
  \bauthor{\bsnm{Nowak},~\bfnm{Adam}\binits{A.}}
(\byear{2002}).
\btitle{Gradient estimates for harmonic and {$q$}-harmonic functions of
  symmetric stable processes}.
\bjournal{Illinois J. Math.}
\bvolume{46}
\bpages{541--556}.
\bid{mr={1936936}}
\end{barticle}
\endbibitem

%b16 ###
\bibitem{MR2013738}
\begin{barticle}[mr]
\bauthor{\bsnm{Bogdan},~\bfnm{Krzysztof}\binits{K.}},
  \bauthor{\bsnm{St{\'o}s},~\bfnm{Andrzej}\binits{A.}} \AND
  \bauthor{\bsnm{Sztonyk},~\bfnm{Pawe{\l}}\binits{P.}}
(\byear{2003}).
\btitle{Harnack inequality for stable processes on {$d$}-sets}.
\bjournal{Studia Math.}
\bvolume{158}
\bpages{163--198}.
\bid{doi={10.4064/sm158-2-5}, mr={2013738}}
\end{barticle}
\endbibitem

%b17 ###
\bibitem{MR2320691}
\begin{barticle}[mr]
\bauthor{\bsnm{Bogdan},~\bfnm{Krzysztof}\binits{K.}} \AND
  \bauthor{\bsnm{Sztonyk},~\bfnm{Pawe{\l}}\binits{P.}}
(\byear{2007}).
\btitle{Estimates of the potential kernel and {H}arnack's inequality for the
  anisotropic fractional {L}aplacian}.
\bjournal{Studia Math.}
\bvolume{181}
\bpages{101--123}.
\bid{doi={10.4064/sm181-2-1}, mr={2320691}}
\end{barticle}
\endbibitem

%b18 ###
\bibitem{MR2256481}
\begin{barticle}[mr]
\bauthor{\bsnm{Bogdan},~\bfnm{K.}\binits{K.}} \AND \bauthor{\bsnm{{\.
  Z}ak},~\bfnm{T.}\binits{T.}}
(\byear{2006}).
\btitle{On {K}elvin transformation}.
\bjournal{J. Theoret. Probab.}
\bvolume{19}
\bpages{89--120}.
\bid{doi={10.1007/s10959-006-0003-8}, mr={2256481}}
\end{barticle}
\endbibitem

%b19 ###
\bibitem{CKS2008}
\begin{bmisc}[auto:SpringerTagBib|2009-01-14|16:51:27]
\bauthor{\bsnm{Chen},~\bfnm{Z.~Q.}\binits{Z.~Q.}},
  \bauthor{\bsnm{Kim},~\bfnm{P.}\binits{P.}} \AND
  \bauthor{\bsnm{Song},~\bfnm{R.}\binits{R.}}
  (\byear{2010}).
Heat kernel estimates for {D}irichlet fractional {L}aplacian. \textit{J. European Math. Soc.}
To appear.
\end{bmisc}
\endbibitem

%b20 ###
\bibitem{MR2357678}
\begin{barticle}[mr]
\bauthor{\bsnm{Chen},~\bfnm{Zhen-Qing}\binits{Z.-Q.}} \AND
  \bauthor{\bsnm{Kumagai},~\bfnm{Takashi}\binits{T.}}
(\byear{2008}).
\btitle{Heat kernel estimates for jump processes of mixed types on metric
  measure spaces}.
\bjournal{Probab. Theory Related Fields}
\bvolume{140}
\bpages{277--317}.
\bid{doi={10.1007/s00440-007-0070-5}, mr={2357678}}
\end{barticle}
\endbibitem

%b21 ###
\bibitem{MR1654824}
\begin{barticle}[mr]
\bauthor{\bsnm{Chen},~\bfnm{Zhen-Qing}\binits{Z.-Q.}} \AND
  \bauthor{\bsnm{Song},~\bfnm{Renming}\binits{R.}}
(\byear{1998}).
\btitle{Estimates on {G}reen functions and {P}oisson kernels for symmetric
  stable processes}.
\bjournal{Math. Ann.}
\bvolume{312}
\bpages{465--501}.
\bid{doi={10.1007/s002080050232}, mr={1654824}}
\end{barticle}
\endbibitem

%b22 ###
\bibitem{2009ZCJT}
\begin{barticle}[mr]
\bauthor{\bsnm{Chen},~\bfnm{Zhen-Qing}\binits{Z.-Q.}} \AND
\bauthor{\bsnm{Tokle},~\bfnm{J.}\binits{J.}}
(\byear{2009}).
\btitle{Global heat kernel estimates for fractional Laplacians in
unbounded open sets}.
\bjournal{Probab. Theory Related Fields}.
DOI:
\href{http://dx.doi.org/10.1007/s00440-009-0256-0}{10.1007/s00440-009-0256-0}.
\bnote{To appear.}
\end{barticle}
\endbibitem

%b23 ###
\bibitem{MR2430977}
\begin{barticle}[mr]
\bauthor{\bsnm{Grigor'yan},~\bfnm{Alexander}\binits{A.}} \AND
  \bauthor{\bsnm{Hu},~\bfnm{Jiaxin}\binits{J.}}
(\byear{2008}).
\btitle{Off-diagonal upper estimates for the heat kernel of the {D}irichlet
  forms on metric spaces}.
\bjournal{Invent. Math.}
\bvolume{174}
\bpages{81--126}.
\bid{doi={10.1007/s00222-008-0135-9}, mr={2430977}}
\end{barticle}
\endbibitem

%b24 ###
\bibitem{MR2417435}
\begin{barticle}[mr]
\bauthor{\bsnm{Grzywny},~\bfnm{Tomasz}\binits{T.}} \AND
  \bauthor{\bsnm{Ryznar},~\bfnm{Micha{\l}}\binits{M.}}
(\byear{2007}).
\btitle{Estimates of {G}reen functions for some perturbations of fractional
  {L}aplacian}.
\bjournal{Illinois J. Math.}
\bvolume{51}
\bpages{1409--1438}.
\bid{mr={2417435}}
\end{barticle}
\endbibitem

%b25 ###
\bibitem{MR2386098}
\begin{barticle}[mr]
\bauthor{\bsnm{Grzywny},~\bfnm{Tomasz}\binits{T.}} \AND
  \bauthor{\bsnm{Ryznar},~\bfnm{Micha{\l}}\binits{M.}}
(\byear{2008}).
\btitle{Two-sided optimal bounds for {G}reen functions of half-spaces for
  relativistic {$\alpha$}-stable process}.
\bjournal{Potential Anal.}
\bvolume{28}
\bpages{201--239}.
\bid{doi={10.1007/s11118-007-9071-3}, mr={2386098}}
\end{barticle}
\endbibitem

%b26 ###
\bibitem{MR2207878}
\begin{barticle}[mr]
\bauthor{\bsnm{Hansen},~\bfnm{Wolfhard}\binits{W.}}
(\byear{2006}).
\btitle{Global comparison of perturbed {G}reen functions}.
\bjournal{Math. Ann.}
\bvolume{334}
\bpages{643--678}.
\bid{doi={10.1007/s00208-005-0719-2}, mr={2207878}}
\end{barticle}
\endbibitem

%b27 ###
\bibitem{MR0142153}
\begin{barticle}[mr]
\bauthor{\bsnm{Ikeda},~\bfnm{Nobuyuki}\binits{N.}} \AND
  \bauthor{\bsnm{Watanabe},~\bfnm{Shinzo}\binits{S.}}
(\byear{1962}).
\btitle{On some relations between the harmonic measure and the {L}\'evy measure
  for a certain class of {M}arkov processes}.
\bjournal{J. Math. Kyoto Univ.}
\bvolume{2}
\bpages{79--95}.
\bid{mr={0142153}}
\end{barticle}
\endbibitem

%b28 ###
\bibitem{MR1991120}
\begin{barticle}[mr]
\bauthor{\bsnm{Jakubowski},~\bfnm{Tomasz}\binits{T.}}
(\byear{2002}).
\btitle{The estimates for the {G}reen function in {L}ipschitz domains for the
  symmetric stable processes}.
\bjournal{Probab. Math. Statist.}
\bvolume{22}
\bpages{419--441}.
\bid{mr={1991120}}
\end{barticle}
\endbibitem

%b29 ###
\bibitem{MR1490808}
\begin{barticle}[mr]
\bauthor{\bsnm{Kulczycki},~\bfnm{Tadeusz}\binits{T.}}
(\byear{1997}).
\btitle{Properties of {G}reen function of symmetric stable processes}.
\bjournal{Probab. Math. Statist.}
\bvolume{17}
\bpages{339--364}.
\bid{mr={1490808}}
\end{barticle}
\endbibitem

%b30 ###
\bibitem{MR1643611}
\begin{barticle}[mr]
\bauthor{\bsnm{Kulczycki},~\bfnm{Tadeusz}\binits{T.}}
(\byear{1998}).
\btitle{Intrinsic ultracontractivity for symmetric stable processes}.
\bjournal{Bull. Polish Acad. Sci. Math.}
\bvolume{46}
\bpages{325--334}.
\bid{mr={1643611}}
\end{barticle}
\endbibitem

%b31 ###
\bibitem{2009KKMS}
\begin{bmisc}[mr]
\bauthor{\bsnm{Kulczycki},~\bfnm{T.}\binits{T.}},
\bauthor{\bsnm{Kwa\'{s}nicki},~\bfnm{M.}\binits{M.}},
  \bauthor{\bsnm{Ma\l ecki},~\bfnm{J.}\binits{J.}} \AND
  \bauthor{\bsnm{St\'{o}s},~\bfnm{A.}\binits{A.}}
(\byear{2009}).
\btitle{Spectral properties
of the Cauchy process on half-line and interval}.
\textit{Proc. London Math. Soc.}
DOI:
\href{http://dx.doi.org/10.1112/plms/pdq010}{10.1112/plms/pdq010}.
To appear.
\end{bmisc}
\endbibitem

%b32 ###
\bibitem{MR2231884}
\begin{barticle}[mr]
\bauthor{\bsnm{Kulczycki},~\bfnm{Tadeusz}\binits{T.}} \AND
  \bauthor{\bsnm{Siudeja},~\bfnm{Bart{\l}omiej}\binits{B.}}
(\byear{2006}).
\btitle{Intrinsic ultracontractivity of the {F}eynman--{K}ac semigroup for
  relativistic stable processes}.
\bjournal{Trans. Amer. Math. Soc.}
\bvolume{358}
\bpages{5025--5057 (electronic)}.
\bid{doi={10.1090/S0002-9947-06-03931-6}, mr={2231884}}
\end{barticle}
\endbibitem

%b33 ###
\bibitem{2009-MK-pa}
\begin{barticle}[mr]
\bauthor{\bsnm{Kwa{\'s}nicki},~\bfnm{Mateusz}\binits{M.}}
(\byear{2009}).
\btitle{Intrinsic ultracontractivity for stable semigroups on unbounded open
  sets}.
\bjournal{Potential Anal.}
\bvolume{31}
\bpages{57--77}.
\bid{doi={10.1007/s11118-009-9125-9}, mr={2507446}}
\end{barticle}
\endbibitem

%b34 ###
\bibitem{luks-2009}
\begin{bmisc}[auto:SpringerTagBib|2009-01-14|16:51:27]
\bauthor{\bsnm{Luks},~\bfnm{T.}\binits{T.}}
(\byear{2009}).
\btitle{Harmonic {H}ardy spaces on smooth domains}. Preprint.
  \bnote{Available at} \url{http://arXiv.org/abs/0909.3370v1}.
\end{bmisc}
\endbibitem

%b35 ###
\bibitem{MR2213639}
\begin{barticle}[mr]
\bauthor{\bsnm{Michalik},~\bfnm{Krzysztof}\binits{K.}}
(\byear{2006}).
\btitle{Sharp estimates of the {G}reen function, the {P}oisson kernel and the
  {M}artin kernel of cones for symmetric stable processes}.
\bjournal{Hiroshima Math. J.}
\bvolume{36}
\bpages{1--21}.
\bid{mr={2213639}}
\end{barticle}
\endbibitem

%b36 ###
\bibitem{MR2238934}
\begin{barticle}[mr]
\bauthor{\bsnm{Rao},~\bfnm{Murali}\binits{M.}},
  \bauthor{\bsnm{Song},~\bfnm{Renming}\binits{R.}} \AND
  \bauthor{\bsnm{Vondra{\v{c}}ek},~\bfnm{Zoran}\binits{Z.}}
(\byear{2006}).
\btitle{Green function estimates and {H}arnack inequality for subordinate
  {B}rownian motions}.
\bjournal{Potential Anal.}
\bvolume{25}
\bpages{1--27}.
\bid{doi={10.1007/s11118-005-9003-z}, mr={2238934}}
\end{barticle}
\endbibitem

%b37 ###
\bibitem{bibRm}
\begin{barticle}[auto:SpringerTagBib|2009-01-14|16:51:27]
\bauthor{\bsnm{Riesz},~\bfnm{M.}\binits{M.}}
(\byear{1938}).
\btitle{Int{\'e}grales de {R}iemann--{L}iouville et potentiels}.
  \bjournal{Acta Sci. Math. Szeged}.
  \bvolume{9}
  \bpages{1--42}.
\end{barticle}


%b38 ###
\bibitem{MR1739520}
\begin{bbook}[mr]
\bauthor{\bsnm{Sato},~\bfnm{Ken-iti}\binits{K.-i.}}
(\byear{1999}).
\btitle{L\'evy Processes and Infinitely Divisible Distributions}.
\bseries{Cambridge Studies in Advanced Mathematics}
\bvolume{68}.
\bpublisher{Cambridge Univ. Press}, \baddress{Cambridge}.
\bid{mr={1739520}}
\end{bbook}
\endbibitem

%b39 ###
\bibitem{MR2255353}
\begin{barticle}[mr]
\bauthor{\bsnm{Siudeja},~\bfnm{Bart{\l}omiej}\binits{B.}}
(\byear{2006}).
\btitle{Symmetric stable processes on unbounded domains}.
\bjournal{Potential Anal.}
\bvolume{25}
\bpages{371--386}.
\bid{doi={10.1007/s11118-006-9022-4}, mr={2255353}}
\end{barticle}
\endbibitem

%b40 ###
\bibitem{MR1719233}
\begin{barticle}[mr]
\bauthor{\bsnm{Song},~\bfnm{Renming}\binits{R.}} \AND
  \bauthor{\bsnm{Wu},~\bfnm{Jang-Mei}\binits{J.-M.}}
(\byear{1999}).
\btitle{Boundary {H}arnack principle for symmetric stable processes}.
\bjournal{J. Funct. Anal.}
\bvolume{168}
\bpages{403--427}.
\bid{doi={10.1006/jfan.1999.3470}, mr={1719233}}
\end{barticle}
\endbibitem

%b41 ###
\bibitem{MR1969798}
\begin{barticle}[mr]
\bauthor{\bsnm{Varopoulos},~\bfnm{N.~Th.}\binits{N.~T.}}
(\byear{2003}).
\btitle{Gaussian estimates in {L}ipschitz domains}.
\bjournal{Canad. J. Math.}
\bvolume{55}
\bpages{401--431}.
\bid{mr={1969798}}
\end{barticle}
\endbibitem

%b42 ###
\bibitem{MR1900329}
\begin{barticle}[mr]
\bauthor{\bsnm{Zhang},~\bfnm{Qi~S.}\binits{Q.~S.}}
(\byear{2002}).
\btitle{The boundary behavior of heat kernels of {D}irichlet {L}aplacians}.
\bjournal{J. Differential Equations}
\bvolume{182}
\bpages{416--430}.
\bid{doi={10.1006/jdeq.2001.4112}, mr={1900329}}
\end{barticle}
\endbibitem

%b43 ###
\bibitem{MR842803}
\begin{barticle}[mr]
\bauthor{\bsnm{Zhao},~\bfnm{Zhong~Xin}\binits{Z.~X.}}
(\byear{1986}).
\btitle{Green function for {S}chr\"odinger operator and conditioned
  {F}eynman--{K}ac gauge}.
\bjournal{J. Math. Anal. Appl.}
\bvolume{116}
\bpages{309--334}.
\bid{doi={10.1016/S0022-247X(86)80001-4}, mr={842803}}
\end{barticle}
\endbibitem

\end{thebibliography}
\end{document}